





\input amstex
\documentstyle{amsppt}
\magnification=\magstep1
\pageheight{9truein}
\pagewidth{6.7truein}
\NoBlackBoxes
\def\ds{\displaystyle}

\def\ep{\varepsilon} 
\def\gA{{\frak A}}
\def\gB{{\frak B}}
\def\rC{\text{\rm C}}
\def\B{{\Cal B}}
\def\D{{\Cal D}}
\def\F{{\Cal F}}
\def\N{{\Cal N}}
\def\nat{{\Bbb N}}
\def\real{{\Bbb R}}
\def\normm{{|\!|\!|}}  
\def\dotnorm{|\cdot|}
\def\dotNorm{\|\cdot\|}
\def\Bar#1{{\skew2\bar{\bar#1}}}
\def\Im{\operatorname{Im}}
\def\supp{\operatorname{supp}}
\def\To{\Rightarrow} 

\topmatter
\title A Banach space block finitely universal for monotone bases\endtitle
\rightheadtext{UNIVERSAL BANACH SPACE}
\author E. Odell and Th. Schlumprecht\endauthor
\thanks Research supported by NSF and TARP.\endthanks 
\abstract A reflexive Banach space $X$ with a basis $(e_i)$ is constructed 
having the property that every monotone basis is block finitely representable 
in each block basis of $X$.\endabstract 
\endtopmatter

\document 
\baselineskip=18pt 

\head \S1. Introduction\endhead 

B.~Maurey and H.~Rosenthal \cite{MR} in 1977 gave an example of a normalized
weakly null sequence for which the summing basis is block finitely 
represented in all subsequences. 
They asked whether one could for some basis achieve this in all block bases 
and thereby solve in the negative the unconditional basic sequence problem. 
The later problem was subsequentially solved in 1991 in a spectacular paper 
of W.T.~Gowers and B.~Maurey \cite{GM}. 
Further examples have since been given solving a myriad of related problems.  
These examples all have a basis $(e_i)$ with the property that for some 
sequence $c_n\nearrow \infty$ if $(y_i)$ is any block basis of $(e_i)$ then 
for all $n\in\nat$ and $\ep>0$ there is a block basis $(z_i)_{i=1}^n$ of 
$(y_i)$ with 
$${\|\sum_{i=1}^n z_i\| \over \|\sum_{i=1}^n (-1)^i z_i\|} >c_n-\ep\ .$$ 
In known examples $c_n$ grows rather slowly (e.g., like $\log n$). 
In this paper we answer in the affirmative the problem of B.~Maurey and 
H.~Rosenthal. 
Thus we obtain an example for which one has $c_n=n$, the worst possible 
unconditionality. 

Our example is a ``conditional version'' of an example in \cite{OS} in which 
an unconditional basis is constructed so that all finite 1-unconditional 
bases are block finitely represented in all block bases. 
The norm in that example satisfies the implicit equation 
$$\|x\| = \max\left\{ \|x\|_\infty,\sup \biggl\{ \dfrac1{f(\ell)} 
\sum_{j=1}^\ell \normm E_i(x)\normm_{m_i} : 
{\ell \in\nat,\ m_i\in\nat,\ \text{for }i=1,\ldots,\ell \atop 
(m_i,E_i)_{i=1}^\ell \text{ is admissible}} \biggr\}\right\}\ ,
\tag 1.1$$
where $\normm x\normm_m \equiv \sup \{\frac1m \sum_{i=1}^m \|F_i(x)\| : 
F_1<\cdots <F_m\}$. 
In (1.1) $f(\ell) \equiv \log_2 (\ell+1)$ and $(m_i,E_i)_{i=1}^\ell$ is 
{\it admissible\/} if $m_1>2$, $f(m_{i+1}) >\sum_{j=1}^i |E_j|$, and 
$E_1< E_2<\cdots <E_\ell$ (the notation will be described in \S2 below). 

The idea of ``defining'' a norm by an implicit equation is due to T.~Figiel 
and W.B.~Johnson \cite{FJ} in their paper on Tsirelson's space. 
We might describe the norm in (1.1) as a ``2-layer'' Tsirelson type norm. 
The example in this paper is a ``3-layer'' norm. 
Also terms like ``$\|F_i(x)\|$'' need to be replaced by ``$x_i^*(x)$'' 
where the $x_i$'s are rather carefully chosen in the manner of \cite{GM}. 
It is also possible to define the norms in (1.1) and below by describing 
the dual ball much in the same manner as Tsirelson did in his original 
paper \cite{T}. 

\head \S2. Terminology and useful facts\endhead 

The linear space of all finitely supported real valued functions on $\nat$ 
is denoted by $c_{00}$. 
For nonempty $E,F\subseteq \nat$, ``$E<F$'' means that $\max E<\min F$. 
$|E|$ is the cardinality of $E$. 
For a sequence of nonzero elements $(x_i) \subseteq c_{00}$, 
``$x_1< x_2<\cdots$'' means that $\supp x_i<\supp x_2<\cdots$ where 
$\supp x_i = \{j\in\nat: x_i(j)\ne0\}$. 
If $x\in c_{00}$ and $E\subseteq \nat$ then $Ex\in c_{00}$ is defined 
by $Ex(j) = x(j)$ if $j\in E$ and $0$ otherwise. 
$(e_i)_1^\infty$ is the unit vector basis for $c_{00}$. 
$(e_i^*)_1^\infty$ are the biorthogonal functionals to $(e_i)$. 

We shall be interested in a collection  $\N$ of certain type norms on $c_{00}$. 
$\dotNorm \in\N$ if $(e_i)$ is a normalized monotone basis for the 
completion of $(c_{00},\dotNorm )$ and $\|e_i^*\|=1$ for all $i$. 
The latter just says that $\|\sum a_i e_i\| \ge \max_i |a_i| \equiv 
\|(a_i)\|_\infty$ for all $(a_i)\in c_{00}$. 
All $\ell_p$ norms and the Tsirelson norm are in $\N$. 
There is a natural partial order on $\N$: 
$\dotNorm \le \dotnorm $ if $\|x\| \le |x|$ for all $x\in c_{00}$. 
Clearly $\dotNorm_{c_0} \le \dotnorm  \le \dotNorm_{\ell_1}$ if 
$\dotnorm  \in\N$. 

If $\dotNorm \in \N$ and $X= (c_{00},\dotNorm )$ then $B(X^*)$, the unit 
ball of $X^*$, is the weak* closure of $\{\sum_1^n a_i e_i^* :n\in\nat$, 
$(a_i)_1^n\subseteq \real\} \cap B(X^*)$ and as such may be identified with 
$$\overline{\Bigl\{ (a_i)\in c_{00} :\|\sum a_i e_i^*\| \le 1\Bigr\} } 
\subseteq [-1,1]^\nat \ ,$$
where $[-1,1]^\nat$ is given the product topology. 

There are a number of Tsirelson type spaces in the literature described 
by implicit equations much like (1.1) (\cite{CS}, \cite{AD}, \cite{GM}, 
\cite{OS}). 
These implicit norms are fixed points of certain mappings on $\N$. 
The next proposition gives a general argument for their existence. 

\proclaim{Proposition 2.1} 
Let $P:\N\to \N$ be order preserving $(\dotnorm\le \dotNorm  \To P \dotnorm  
\le P\dotNorm )$. 
Then $P$ admits a smallest fixed point. 
Thus there exists $\dotNorm \in \N$ with 
\roster 
\item"a)" $\dotNorm = P\dotNorm $ 
\item"b)" If $P\dotnorm  =\dotnorm $ then $\dotNorm \le \dotnorm $. 
\endroster
\endproclaim 

\demo{Proof}
Let $\dotNorm_0 \equiv \dotNorm_\infty$. 
By transfinite induction we define $\dotNorm_\alpha$ for $\alpha <\omega_1$. 
If $\alpha = \beta+1$ we set 
$$\|x\|_\alpha = \max (\|x\|_\beta, P\|x\|_\beta)\ .$$ 
If $\alpha$ is a limit ordinal we set $\|x\|_\alpha = \sup_{\beta<\alpha} 
\|x\|_\beta$. 
Clearly $\dotNorm_\alpha \in \N$ for all $\alpha <\omega_1$. 
Since $(\dotNorm_\alpha)_{\alpha<\omega_1}$ is an increasing family of 
norms the dual balls $B_\alpha^* \equiv B(c_{00},\dotNorm_\alpha^*)$ are 
increasing closed subsets of $[-1,1]^\nat$. 
Since this space is compact metrizable there exists $\alpha_0<\omega_1$ 
so that $B_{\alpha_0}^* = B_\beta^*$ for all $\beta \ge \alpha_0$. 
Thus $\dotNorm_{\alpha_0} = \dotNorm_\beta$ for $\beta \ge\alpha_0$ as well.  
We set $\dotNorm= \dotNorm_{\alpha_0}$. 

To see that $\dotNorm$ is a fixed point for $P$ we first note that 
$\dotNorm = \dotNorm_{\alpha_0+1} \ge P\dotNorm_{\alpha_0} = P\dotNorm$. 
For the reverse inequality it suffices to show by induction that 
$\dotNorm_\alpha \le P\dotNorm$ for all $\alpha <\omega_1$. 
Clearly this holds for $\alpha=0$ and if $\dotNorm_\beta \le P\dotNorm$ 
then $\dotNorm_{\beta+1} = \max (\dotNorm_\beta,P\dotNorm_\beta) \le 
P\dotNorm$. 
Indeed $\dotNorm_\beta \le \dotNorm$ by the definition of $\dotNorm$ and so 
$P\dotNorm_\beta \le P\dotNorm$ since $P$ is order preserving. 
Also if $\alpha$ is a limit ordinal and $\dotNorm_\beta \le P\dotNorm$ 
for $\beta <\alpha$ then $\dotNorm_\alpha \le P\dotNorm$. 
Thus $P\dotNorm = \dotNorm$. 

To see b), let $P\dotnorm =\dotnorm$. 
Then $\dotNorm_0 \le \dotnorm$ and by induction we easily obtain 
$\dotNorm_\alpha \le\dotnorm$ for $\alpha<\omega_1$, hence 
$\dotNorm \le \dotnorm$.\qed
\enddemo

Basic sequences $(x_i)$ and $(y_i)$ are {\it $C$-equivalent\/} if for some 
constants $\alpha,\beta$ with $\alpha^{-1}\beta \le C$, 
$\alpha \| \sum a_ix_i\| \le \|\sum a_i y_i\| \le \beta \|\sum a_ix_i\|$ 
for all $(a_i) \in c_{00}$. 
A basic sequence $(z_i)$ is {\it block finitely represented\/} in a basic 
sequence $(x_i)$ if for all $\ep>0$ and $n\in\nat$ there exists a block basis 
$(y_i)_1^n$ of $(x_i)$ which is $(1+\ep)$-equivalent to $(z_i)_1^n$. 

Let $\dotNorm \in \N$ and $X= (c_{00},\dotNorm)$. 
For $A\subseteq X^*$ and $x\in X$ we set $\|x\|_A=\sup\{|x^*(x)| :x^*\in A\}$. 
For $1\le p\le\infty$, $C\ge1$ and $k\in\nat$, $x\in X$ is called an 
{\it $\ell_p^k$-average with constant\/} $C$ if $x= k^{-1/p} \sum_{i=1}^k x_i$ 
for some normalized sequence $x_1<\cdots < x_k$ which is $C$-equivalent to 
the unit vector basis of $\ell_p^k$. 

If $\dotNorm \in \N$ and $X= (c_{00},\dotNorm)$ then $\dotNorm^* \in \N$ 
as well. 
Indeed $(e_i^*)$ is a normalized monotone basic sequence in $X^*$ and 
$\|\sum a_i e_i^*\| \ge \max |a_i|$ for $(a_i) \in c_{00}$. 
So we are free to use our notation $x^* <y^*$ for elements of the dual 
in $\text{span}(e_i^*)$ as well. 

Before defining our norm we present some technical notation and a lemma. 
The lemma could be postponed but it does help one become familiar with 
the terminology. 
Our first definitions are motivated by \cite{GM}. 

Fix $H\subset c_{00} \cap [-1,1]^\nat$, a countable subset of nonzero elements. 
Let $\D= (D_i)$ be a sequence of subsets of $H$. 
Let 
$$\sigma :\{ (x_1^*,\ldots,x_n^*) : n\in\nat\ ,
\ x_i^*\in H\text{ for } i\le n\ , x_1^*<\cdots <x_n^*\} \to \nat$$ 
be an injective map satisfying the following condition 
$$\sigma(x_1^*,\ldots,x_n^*) >\max (k,\max\supp x_n^*) 
\text{ if }x_n^* \in D_k\ .
\tag $\sigma,\D$ $$ 
Let $M= (M_n)$ be a subsequence of $\nat$. 

\definition{Definition} 
$(x_1^*,\ldots,x_n^*)\subseteq c_{00}$ is {\it $(\D,M,\sigma)$-admissible\/} if 
\roster
\item"1)" $x_1^* <\cdots < x_n^*$
\item"2)" $x_1^* \in \bigcup_{i\ge M_n} D_i$ 
\item"3)" $x_{i+1}^* \in D_{\sigma (x_1^*,\ldots,x_i^*)}$  if $1\le i<n$. 
\endroster
\enddefinition 

We have used ``$x_i^*$'' above in our definitions because we will be applying 
this for elements in $X^*$. 

Let $f:[1,\infty)\to\real$ be given by $f(x) = \log_2(x+1)$. 
We will make use of the fact that $f$ is strictly increasing, $f(1)=1$ and 
both $f(x)$ and $x\over f(x)$ are concave functions. 

For $n\in\nat$ set 
$$\Gamma_n=\Gamma_n (\D,M,\sigma) = \biggl\{ {1\over f(n)} \sum_{j=1}^n : 
(x_j^*)_1^n\text{ is } (\D,M,\sigma)\text{-admissible}\biggr\}\ .$$ 
Let $\Gamma = \Gamma (\D,M,\sigma) = \bigcup_\nat \Gamma_n$, 
$D= \bigcup_{i\ge M_1} D_i$ and note that $\Gamma_1=D$. 

\proclaim{Lemma 2.2} 
Let $\dotNorm \in\N$, $X= (c_{00},\dotNorm)$ and suppose that $\Gamma (\D,M,
\sigma) \subseteq B(X^*)$. 

{\rm a)} Let $k,m\in\nat$, $\ep>0$ and let $e_k<y_1<\cdots <y_m$ be a block 
sequence of nonzero elements in $B(X) \cap [-1,1]^\nat$. 
Assume that for all $1<i\le m$ and for any $(\D,M,\sigma)$-admissible sequence 
$(x_1^*,\ldots,x_j^*)$ with $\max\supp x_j^*\ge \min \supp y_i$ we have 
for all $x^*\in \bigcup_{i\ge \sigma (x_1^*,\ldots,x_j^*)} D_i$,   
$$|x^* (y_i)| \le {\ep \over f^{-1}(\frac1{\ep} \max\supp y_{i-1})}\ .
\tag 2.2.1$$

Let $\ell\in\nat$, $x^* = \frac1{f(\ell)} \sum_{j=1}^\ell x_j^* \in\Gamma_\ell$ 
and $y=\sum_{i=1}^m \alpha_i y_i$ with $(\alpha_i)_1^m \subseteq \real$ so 
that $x^*(y)\ne0$. 
Set 
$$j_1 \equiv \min \{ j\le \ell :\max\supp x_j^* \ge\min\supp y\}\ .$$
Then
$$\leqalignno{\qquad	
|x^*(y)| & \le {1\over f(\ell)} |x_{j_1}^* (y)| +\max_{i\le m} |\alpha_i|\ 
\Big| \left( x^* - {x_{j_1}^*\over f(\ell)} \right) y_i\Big|   &(2.2.2)\cr 
&\qquad 
+ {\min (m-1,\ell-1)\over f(\ell)} \sup 
\biggl\{ |x^*(y)| : x^* \in \bigcup_{t\ge k} D_t\biggr\}\cr 
&\qquad 
+ {\min (2,\ell-1)\over f(\ell)} \sup 
\biggl\{ |x^*(\alpha_iy_i)|: i\le m\ ,\ x^*\in \bigcup_{t\ge k} D_t\biggr\}
+ 2\ep \| (\alpha_i)\|_\infty \cr 
&\le {1\over f(\ell)} |x_{j_1}^* (y)|  
+ {\min (m-1,\ell-1)\over f(\ell)} 
\sup \biggl\{ |x^*(y)| : x^* \in \bigcup_{t\ge k} D_t\biggr\}\cr 
&\qquad 
+ \|(\alpha_i)\|_\infty \Biggl[ \min (\ell-1,1) \left( 1+{1\over f(\ell)}\right)
	 +2\ep \cr 
&\hskip2.0truein + \sup \biggl\{ |x^* (y_i)| : i\le m\ ,\ 
	x^* \in \bigcup_{t\ge k} D_t\biggr\} \Biggr]
\cr}$$

{\rm b)} Let $(y_i)_{i=1}^\infty \subseteq B(X) \cap [-1,1]^\nat$ be a 
block basis of $(e_i)$ satisfying 
$$\ep_n \equiv \sup \biggl\{ |x^* (y_i)| : i\ge n\ ,\ x^* \in \bigcup_{t\ge n} 
D_t\biggr\}\to 0
\tag 2.2.3$$ 
as $n\to\infty$.

Then  for all $\ell,m\in \nat$ and $(\alpha_i)_1^m \subseteq \real$ we have 
$$\leqalignno{ \qquad
&\varlimsup_{n_1\to\infty} \ldots \varlimsup_{n_m\to\infty} 
\Big\| \sum_{i=1}^m \alpha_i y_{n_i} \Big\|_{\Gamma_\ell} 
\le {1\over f(\ell)} \varlimsup_{n_1\to\infty}\ldots \varlimsup_{n_m\to\infty} 
	\Big\|\sum_{i=1}^m \alpha_i y_{n_i}\Big\|_D\ \text{ and}
&(2.2.4)\cr 
&\varlimsup_{n_1\to\infty} \ldots\varlimsup_{n_m\to\infty} 
	\Big\|\sum_{i=1}^m \alpha_i y_{n_i}\Big\|_\Gamma 
\le \max \biggl\{ \|(\alpha_i)\|_\infty , \varlimsup_{n_1\to\infty}\ldots 
\varlimsup_{n_m\to\infty} \Big\|\sum_{i=1}^m \alpha_i y_{n_i}\Big\|_D
\biggr\}\ .
&(2.2.5)\cr}$$
\endproclaim

\demo{Proof} 
a) We begin by choosing for $1\le i\le m$, $\tilde I_i \subseteq \{j_1+1,
\ldots ,\ell\}$ to be those $j$'s for which  $x_j^*$ acts only on $y_i$. 
For $1<i<m$ 
$$\tilde I_i = \{j>j_1 :\supp x_j^* \cap \supp y_i\ne\emptyset\text{ and } 
y_{i-1} <x_j^* <y_{i+1}\}\ .$$ 
We take 
$$\align
\tilde I_1 & = \{j>j_1 :\supp x_j^* \cap \supp y_1\ne\emptyset\text{ and } 
x_j^* <y_2\}\ ,\cr
\tilde I_m & = \{j>j_1: \supp x_j^*\cap \supp y_m\ne\emptyset\text{ and } 
y_{m-1} <x_j^*\}\ .\endalign$$ 
For $2\le i\le m$ it may happen that there exists (at most one) $j_i>j_1$ 
with $\min\supp x_{j_i}^* \le \max\supp y_{i-1}$ and 
$\min\supp y_i \le \max\supp x_{j_i}^*$ and if $i<m$, 
$\max\supp x_{j_i}^* < \min \supp y_{i+1}$. 
In this case we take $I_i=\tilde I_i$. 
Otherwise if $\tilde I_i \ne\emptyset$ we set $j_i = \min\Tilde I_i$ and 
$I_i = \tilde I_i\setminus \{j_i\}$. 
Let $I_0$ be the set of all $j_i$'s, $i\ge2$,  
thus obtained and note that $|I_0|\le \min (\ell-1,m-1)$. 
Let $I_1 = \tilde I_1$. 
We thus have 
$$\{j_1\} \cup \bigcup_{i=0}^m I_i\supseteq \{ j\le \ell :x_j^* (y)\ne0\}\ .$$ 
Finally we set 
$$i_0 = \min \left\{ i\le m:f(\ell) \le {\max\supp y_i\over\ep}\right\}$$ 
and $i_0=m$ if the set is empty. 
Since $y,x^*\in [-1,1]^\nat$, 
$$\sum_{i=1}^{i_0-1} {1\over f(\ell)} \sum_{j\in I_i} |x_j^* (y)| 
\le {\max\supp y_{i_0-1} \over f(\ell)} <\ep \ .
\tag 1$$ 
For $m\ge i>i_0$ and $j\in I_i$ it follows from (2.2.1) and the 
fact that $(x_1^*,\ldots,x_{j-1}^*)$ is $(\D,M,\sigma)$-admissible with 
$\max\supp x_{j-1}^* \ge \min \supp y_i$ that 
$$ |x_j^* (y_i)|  < {\ep \over f^{-1}\left({\max\supp y_{i-1} 
\over \ep}\right)} 
\le {\ep \over f^{-1} \left( {\max\supp y_{i_0} \over \ep} \right)} 
<{\ep\over\ell}$$
where the last inequality uses the definition of $i_0$. 
Thus 
$$\sum_{i=i_0+1}^m {1\over f(\ell)} \sum_{j\in I_i} |x_j^* (y_i)| 
< \sum_{i=i_0+1}^m {1\over f(\ell)} |I_i| {\ep\over\ell} \le\ep\ .
\tag 2$$ 

We are now ready to estimate $|x^* (y)|$. 
$$\eqalign{ |x^*(y)| & = \Big| {x_{j_1}^*(y)\over f(\ell)} +\alpha_{i_0} 
\left( x^* - {x_{j_1}^* \over f(\ell)} \right) (y_{i_0})  
+ {1\over f(\ell)} \sum_{j\in I_0} x_j^* (y-\alpha_{i_0} y_{i_0}) 
+ {1\over f(\ell)} \sum_{i\ne i_0} \sum_{j\in I_i} x_j^* 
(\alpha_i y_i)\Big| \cr 
&\le {1\over f(\ell)} |x_{j_1}^* (y)| + |\alpha_{i_0}| \ 
\left| \left( x^* - {x_{j_1}^* \over f(\ell)} \right) (y_{i_0})\right|
+ {1\over f(\ell)} \sum_{j\in I_0} |x_j^* (y)|\cr  
&\qquad +{1\over f(\ell)} \sum_{j\in I_0} |\alpha_{i_0}|\ |x_j^* (y_{i_0})|
+ 2\ep \|(\alpha_i)\|_\infty\ ,
\cr}$$ 
where the very last estimate follows from (1) and (2). 
Thus 
$$\align 
|x^* (y)| & \le {1\over f(\ell)} |x_{j_1}^* (y)| + \max_{i\le m} 
|\alpha_i| \Big| \left( x^* - {x_{j_1}^* \over f(\ell)}\right)(y_i)\Big|
\tag 3\cr 
&\qquad + {\min (m-1,\ell-1)\over f(\ell)} \sup \biggl\{ |x^* (y)| : 
x^* \in \bigcup_{t\ge k} D_t\biggr\} \cr 
&\qquad +{\min (2,\ell-1)\over f(\ell)} \supp 
\biggl\{ |x^* (\alpha_i y_i)| :i\le m\ ,\ x^* \in \bigcup_{t\ge k} D_i\biggr\}
+ 2\ep \|(\alpha_i)\|_\infty\ .
\endalign$$
We have used that if $j\in I_0$ then $j>j_1$ and so from the fact that 
$(x_1^*,\ldots,x_\ell^*)$ is $(\D,M,\sigma)$-admissible and 
$\max\supp x_{j_1}^*\ge \min \supp y>k$ we have that $x_j^*\in D_t$ for 
some $t>k$. 
Furthermore the ``$\min (2,\ell-1)$'' came from 
$$| \{j\in I_0 : x_j^* (y_{i_0}) \ne0\}| \le \min (2,\ell-1)\ .$$ 
Continuing, the right hand expression in (3) is 
$$\align 
&\le {|x_{j_1}^* (y)|\over f(\ell)} + \min (\ell-1,1) 
\left( 1+{1\over f(\ell)}\right) \|(\alpha_i)\|_\infty \tag 4\cr 
&\qquad + {1\over f(\ell)} \min (m-1,\ell-1) \sup \biggl\{ |x^*(y)| : 
x^* \in \bigcup_{t\ge k} D_t\biggr\} \cr 
&\qquad + {1\over f(\ell)} \min (2,\ell-1) \|(\alpha_i)\|_\infty 
\sup \biggl\{ |x^* (y)| :i\le m\ ,\ x^* \in \bigcup_{t\ge k} D_t\biggr\}\cr 
&\qquad + 2\ep \|(\alpha_i)\|_\infty\ .
\endalign$$
(3) and (4) combined yield (2.2.2)  (using ${\min (2,\ell-1)\over f(\ell)} 
\le 1$). 

b) To deduce (2.2.4) it suffices to prove that given $\ep>0$ and $m\in\nat$ 
$$\exists\ \tilde n_1 \in \nat\ \forall\ n_1\ge \tilde n_1 \ \exists\ 
\tilde n_2 >n_1\ \forall\ n_2 \ge \tilde n_2 \ldots\ \exists\ \tilde n_m 
> n_{m-1} \ \forall\ n_m\ge \tilde n_m$$
so that we have 
\roster
\item"$(*)$" 
for all $\ell \in \nat$ and $x^* = {1\over f(\ell)} 
\sum_{j=1}^\ell x_j^* \in\Gamma_\ell$ and $(\alpha_i)_1^m\subseteq \real$, 
$$\Big| x^* \biggl( \sum_{i=1}^m \alpha_i y_{n_i}\biggr) \Big| 
\le {1\over f(\ell)} \Big\| \sum_{i=1}^m \alpha_i y_{n_i}\Big\|_D 
+ (\ell+4) \ep \|(\alpha_i)\|_\infty\ \text{ and}$$
\item"$(**)$" 
for any  $(\alpha_i)_1^m \subseteq \real$ 
$$\Big\| \sum_{i=1}^m \alpha_i y_{n_i}\Big\|_\Gamma 
\le \max \biggl\{ \|(\alpha_i)\|_\infty ,\Big\| \sum_{i=1}^m \alpha_i 
y_{n_i}\Big\|_D\biggr\} 
+ 7\ep \|(\alpha_i)\|_\infty\ .$$
\endroster

To see this given $m$ and $\ep$ begin by choosing $\tilde n_1$ so that 
$$\align 
&\ep_{\tilde n_1} < {\ep \over 2m^2} \ \text{ and}\tag 5\cr 
&\text{If }\ell\in \nat\text{ with } {\ell\over f(\ell)} \ep_{\tilde n_1} 
> \ep \text{ then } {1\over f(\ell)} < {\ep\over m}\ .\tag 6
\endalign$$
Assume $\tilde n_1 <n_1 <\cdots < \tilde n_i <n_i$ are chosen 
$(n_j >\tilde n_j$ arbitrarily and $\tilde n_j$ dependent on $n_{j-1}$). 
Use (2.2.3) to choose $\tilde n_{i+1} >n_i$ so that 

\noindent\hang (7) 
for any $n_{i+1} \ge \tilde n_{i+1}$ and any $(\D,M,\sigma)$-admissible 
sequence $(x_1^* ,\ldots, x_j^*)$ with $\max\supp x_j^* \ge \min\supp 
y_{n_{i+1}}$ and $x^* \in D_t$ for some $t\ge \sigma (x_1^*, \ldots,x_j^*)$ 
we have 
$$|x^* (y_{n_{i+1}})| < {\ep\over  f^{-1} \left( {\max\supp y_{n_i} \over\ep}
\right)}\ .$$ 

Now let $\ell\in\nat$, $y= \sum_1^m \alpha_i y_{n_i}$ and 
$x^* = {1\over f(\ell)} \sum_{j=1}^\ell x_j^* \in\Gamma_\ell$. 
Let $j_1$ be defined as in a). 
>From the first inequality in (2.2.2) we have ($k$ is replaced by 
$\tilde n_1$) 
$$\align 
|x^* (y)| & \le {1\over f(\ell)} |x_{j_1}^* (y)| +\|(\alpha_i)\|_\infty 
\max_{i\le m} \sum_{j=j_1+1}^\ell {|x_j^* (y_{n_i})| \over f(\ell)}
+ {m\over f(\ell)} \cdot m\cdot \ep_{\tilde n_1} \cdot 
\|(\alpha_i)\|_\infty 
\tag 8\cr 
&\qquad + {2\over f(\ell)} \ep_{\tilde n_1} \|(\alpha_i)\|_\infty 
+ 2\ep \|(\alpha_i)\|_\infty\cr  
&\le {1\over f(\ell)} |x_{j_1}^* (y)| + \|(\alpha_i)\|_\infty \ell 
\ep_{\tilde n_1} 
+ \|(\alpha_i)\|_\infty \cdot 4\ep\endalign$$ 
where the last inequality uses (5) to deduce $m^2 \ep_{\tilde n_1}<\ep$ 
and the trivial estimate ${2\ep_{\tilde n_1}\over f(\ell)} <\ep$. 
Thus
$$|x^*(y)| \le {1\over f(\ell)} \|y\|_D + \|(\alpha_i)\|_\infty 
\ep (\ell+4)\ .$$ 
This completes the proof of $(*)$. 

To see $(**)$ let $\ell\in\nat$ and $x^* = {1\over f(\ell)} \sum_{j=1}^\ell 
x_j^* \in\Gamma_\ell$ so that for $y= \sum_1^m \alpha_i y_{n_i}$, 
$$\|y\| = \|y\|_{\Gamma_\ell} \le |x^* (y)| + \ep \|(\alpha_i)\|_\infty\ .$$

\noindent {\it Case 1\/}. 
${\ell \over f(\ell)} \ep_{\tilde n_1} \le \ep$

>From (8) above we have 
$$\align 
|x^* (y)| & \le {1\over f(\ell)} \|y\|_D + \|(\alpha_i)\|_\infty 
\max_{i\le m} \sum_{j=j_1+1}^\ell {|x_j^* (y_{n_i})|\over f(\ell)} 
+ \|(\alpha_i)\|_\infty \cdot 4\ep\cr 
&\le {1\over f(\ell)} \|y\|_D + \|(\alpha_i)\|_\infty 
\left[ {\ell \ep_{\tilde n_1} \over f(\ell)} + 4\ep\right]\cr 
&\le {1\over f(\ell)} \|y\|_D + \|(\alpha_i)\|_\infty \cdot 5\ep 
\endalign$$ 

\noindent {\it Case 2\/}.  
${\ell\over f(\ell)} \ep_{\tilde n_1} >\ep$ and hence by (6), 
${1\over f(\ell)} < {\ep\over m}$ 

Note that 
$$\align 
\max_{i\le m} |\alpha_i| \Big| \left( x^* - {x_{j_1}^* \over f(\ell)}\right) 
(y_i)\Big| & \le \|(\alpha_i)\|_\infty \left( 1+ \max_{i\le m} 
{|x_{j_1}^* (y_i)|\over f(\ell)} \right)\cr 
&\le \| (\alpha_i)\|_\infty \left( 1+{1\over f(\ell)}\right)\ .
\endalign$$
Thus using (2.2.2) as in (8) except for this estimate we have 
$$\align 
|x^* (y)| & \le {\|y\|_D\over f(\ell)} + \|(\alpha_i)\|_\infty 
\left( 1+{1\over f(\ell)}\right) + \|(\alpha_i)\|_\infty \cdot 4\ep\cr 
&\le \left( {m\over f(\ell)} + 1+{1\over f(\ell)} +4\ep\right) 
\|(\alpha_i)\|_\infty \cr 
&< (6\ep +1) \|(\alpha_i)\|_\infty\ .\endalign$$ 
We used $\|y\|_D \le m\|(\alpha_i)\|_\infty$. 
This completes the proof of $(**)$.\qed
\enddemo

\remark{Remark 2.3} 
Lemma 2.2 also holds in the setting where $\Gamma\equiv \bigcup_{n\ge2}
\Gamma_n$ and we shall also use it this way below. 
\endremark 

\head \S3. The construction of $X$\endhead 

Let $H\subset c_{00}\cap [-1,1]^\nat$ be a countable set of nonzero elements 
satisfying the following three properties. 

\roster
\item"(H$_1$)" $H$ is dense in $c_{00}$ with respect to $\dotNorm_{\ell_1}$ 
\item"(H$_2$)" If $a\in H$, $I$ is an interval of integers and 
$Ia\ne 0$ then $Ia\in H$. 
\item"(H$_3$)" For all $n\in\nat$, $a_1<\cdots < a_n$ in $H$ we have 
$\sum_{i=1}^n a_i$, ${1\over f(n)} \sum_{i=1}^n a_i$ and 
${1\over n} \sum_{i=1}^n a_i$ are all in $H$.
\endroster 

Let $M= (M_n)\subseteq \nat$ be strictly increasing with $M_1=2$. 
Let 
$$\sigma : \{(a_1,\ldots,a_n) : n\in\nat\ ,\ a_1<\cdots < a_n\ ,\ 
a_i \in H\text{ for } i\le n\} \to \nat$$ 
be an injection satisfying the following four properties. 
\roster
\item"$(\sigma_1)$" 
Let $n\in\nat$ and $a_1<\cdots <a_n$ be in $H$. 
Let $I$ be an interval in $\nat$ and suppose $[j_1,j_2] = \{i:Ia_i\ne 0$, 
$i\le n\}\ne\emptyset$. 
Then $\sigma (Ia_{j_1},\ldots,Ia_{j_2}) \le \sigma (a_1,\ldots,a_n)$.
\item"$(\sigma_2)$" 
For $n\in\nat$ and $a_1<\cdots <a_n$ in $H$, 
$$\max\supp a_n <\sigma (a_1,\ldots,a_n)$$ 
\item"$(\sigma_3)$" 
Let $I_m(\sigma)$ be the range of $\sigma$. 
Then 
$$C_f(\sigma) \equiv \sum_{n\in I_m(\sigma)} {1\over f(n)} <\infty\ .$$ 
In particular for $1\le p<\infty$ 
$$C_p (\sigma) \equiv \sum_{n\in I_m (\sigma)} n^{-1/p} <\infty\ .$$ 
\item"$(\sigma_4)$" 
For $m\in \Im (\sigma)$ let $\bar m$ and $\Bar m$ be the predecessor and 
successor, respectively, of $m$ in $\Im (\sigma)$. 
Then if $\ell \in [1,\bar m]\cup [\Bar m,\infty)$ 
$$ 2\left[ {f(m)\over m} + {1\over f(\ell)} + {f(m)\over m} \ 
{\min (\ell,m)\over f(\ell)}\right] 
\le \left\{ \eqalign{&{3\over f(\ell)} \ ,\quad \ell<m\cr 
	&{3f(m)\over m}\ ,\quad \ell>m\cr}
	\right.$$
\endroster

We have made no attempt to give a minimal list of necessary conditions but 
rather have chosen to list precisely the conditions we will use. 
It is straightforward to prove that such a set $H$ and function $\sigma$ 
exist. 

For $\dotNorm \in \N$ let $X= (c_{00},\dotNorm )$. 
For $m\ge 2$ let 
$$A_m^X  = \biggl\{ {1\over m} \sum_{i=1}^m a_i^* :a_i^* \in H 
\cap B(X^*) \
\text{for } i\le m \text{ and } a_1^* <\cdots < a_m^*\biggr\}\ .
$$ 
Set $A^X = \bigcup_{m\ge2} A_m^X$, $\tilde A_m^X = \bigcup_{n\ge m} A_n^X$ 
for $m\ge2$ and $\tilde{\gA}^X = (\tilde A_m^X)_{m=2}^\infty$. 
Note that $A^X \subseteq  B(X^*)$. 
For $m\ge2$ let 
$$B_m^X  = \biggl\{ {1\over f(m)} \sum_{i=1}^m a_i^* : (a_1^*,
\ldots,a_m^*) \subseteq A^X 
\text{is $(\tilde\gA^X,M,\sigma)$-admissible}\biggr\}\ .
$$
Let $B^X = \bigcup_{n=2}^\infty B_n^X$ and $\gB^X = (B_n^X)_{n=2}^\infty$. 
For $m\ge1$ set 
$$C_m^X = \biggl\{ {1\over f(m)} \sum_{i=1}^m a_i^* : (a_1^*,
\ldots, a_m^*) \subseteq B^X
\text{is $(\gB^X,M,\sigma)$-admissible}\biggr\}\ .
$$ 
Note that $C_1^X = B^X$ and let $C^X = \bigcup_{n=1}^\infty C_n^X$. 

In the notation of \S2 
$$B_n^X = \Gamma_n (\tilde\gA^X,M,\sigma)\ \text{ and }\ 
C_n^X = \Gamma_n (\gB^X ,M,\sigma)\ .$$ 

\remark{Remark} 
Any element of $A_n^X$, $\tilde A_n^X$ or $B_n^X$ has at least $n$ 
nonzero coordinates. 
Thus this plus $(\sigma_2)$ implies that conditions $(\sigma,\tilde\gA)$ 
and $(\sigma,\gB)$ hold (see \S2). 
Also $A_n^X,B_n^X, C_n^X \subseteq [-1,1]^\nat$ for all $n$.
\endremark 

\proclaim{Proposition 3.1} 
There exists a norm $\dotNorm \in\N$ so that for $X= (c_{00},\dotNorm)$ 
and all $x\in c_{00}$, 
$$\|x\| = \max (\|x\|_\infty, \|x\|_{C^X})\ .
\tag 3.1.1$$ 
Moreover the completion of $X$ is reflexive.
\endproclaim 

\demo{Proof} 
Define $P:\N\to\N$ by $P\|x\| = \max \{\|x\|_\infty, \|x\|_{C^{X,\dotNorm}}\}$. 
To see that $P\dotNorm \in\N$ we need only note that $(e_i)$ is monotone 
for $P\dotNorm$. 
This follows from the fact that if $I\subseteq \nat$ is an initial 
interval ($I=[1,n]$ for some $n$) and $x^*\in C_\ell^X$ (or $A_\ell^X,
B_\ell^X$) then for all $m>n$ there exists $y^* \in C_\ell^X$ with 
$Iy^* = Ix^*$ and $y^* (k)=0$ if $n<k\le m$. 
This easily established fact uses that $(e_i)$ is a monotone basis 
for $\dotNorm$. 

$P$ is order preserving. 
This follows from the fact that if $\dotnorm \le\dotNorm$ are norms in $\N$ 
then $C^{X,\dotnorm} \subseteq C^{X,\dotNorm}$. 
Thus by Proposition~2.1 we obtain a norm $\dotNorm\in\N$ satisfying (3.1.1). 

It remains to show that the completion of $X=(c_{00},\dotNorm)$ 
is reflexive. 
To do this we shall prove that $(e_i)$ is shrinking and boundedly complete. 
If $(e_i)$ were not boundedly complete then there exists a block basis 
$(x_i)$ of $(e_i)$ with $\|x_i\|\ge1$ for all $i$ and $\|\sum_{i=1}^n x_i\| 
\le K$ for all $n$ and some $K<\infty$. 
Choose $x_i^* \in B(X^*)\cap H$ with for all $i$, $\supp x_i^*\subseteq 
[\min \supp x_i,\max\supp x_i]$ and $x_i^* (x_i) >1/3$. 
Indeed using that $(e_i)$ is a monotone basis for $X$ we can produce $y_i^* 
\in X^*$ with $\|y_i^*\| \le 2$, 
$\supp y_i^* \subseteq [\min\supp x_i,\max\supp x_i]$ and 
$y_i^* (x_i) = \|x_i\| \ge1$. 
Then we take $x_i^*$ in $H\cap B(X^*)$ to be an appropriate approximation 
of $y_i^*\over \|y_i^*\|$. 
Let $\ell>1$ be fixed. 
Then $a_1^* \equiv {1\over M_\ell} \sum_1^{M_\ell} x_i^* \in A_{M_\ell}^X$. 
We choose $a_2^* = {1\over\sigma (a_1^*)} \sum_{i=M_\ell+1}^{M_\ell +\sigma 
(a_1^*)} x_i^*$ and so forth obtaining ultimately $(a_1^*,\ldots,a_\ell^*)$ 
which is $\tilde\gA,M,\sigma)$-admissible. 
It follows that ${1\over f(\ell)} \sum_1^\ell a_\ell^* \in B(X^*)$ 
and for an appropriate $n$ 
$${1\over f(\ell)} \sum_1^\ell a_i^* \left(\sum_{j=1}^n x_i\right) 
\ge {1\over f(\ell)} {\ell\over3}\ .$$ 
But this cannot be always bounded by $K$. 

Finally we prove that $(e_i)$ is shrinking. 
If not there exists a normalized block basis $(x_i)$ of $(e_i)$ and 
$\delta>0$ so that $\|\sum_i a_ix_i\| >\delta \sum_i a_i$ for all 
$(a_i) \in c_{00}^+$ (i.e., $a_i\ge 0$ for all $i$). 
Using James' argument \cite{J} that $\ell_1$ is not distortable for this 
$\ell_1^+$-basis we can, by replacing $(x_i)$ by a normalized block basis, 
assume that for a given $\ep>0$,  
$\|\sum a_ix_i\| >(1-\ep) \sum a_i$ if $(a_i) \in c_{00}^+$. 
If $\tilde C$ is the weak* closure of $C^X \cup \{e_n^*\}_\nat$ then we may 
regard $X\subseteq \rC(\tilde C)$, the space of continuous functions on 
$\tilde C$. 
Let $\F = \{F\subseteq \nat:$ there exists $x^*\in  \tilde C$ with 
$x^* (x_i)>1-2\ep$ for all $i\in F\}$. 
$\F$ is a hereditary family of subsets of $\nat$ and has the additional  
property that if $(a_i) \subseteq c_{00}^+$ with $\sum_i a_i=1$ then there 
exists $F\in\F$ with $\sum_F a_i\ge \ep$. 
[Indeed let $x^* \in \tilde C$ with $x^* (\sum a_ix_i)>1-\ep$ and set 
$F= \{i:x^*(x_i) >1-2\ep\}$. 
Then 
$$1-\ep < x^*\Bigl( \sum a_i x_i\Bigr) \le (1-2\ep) \sum_{i\notin F} a_i 
+ \sum_{i\in F} a_i \le (1-2\ep) + \sum_{i\in F} a_i$$ 
and so $\ep \le \sum_{i\in F} a_i$.] 

Ptak's theorem (\cite{P}, see also \cite{BHO}) 
yields that there exists a subsequence $N$ of 
$\nat$ so that $F\in \F$ for all $F\subseteq N$. 
Thus relabeling $(x_i)_N$ as $(x_i)$ we have that for all $n$ there exists 
$x_n^* \in C^X\cup \{e_m^*\}_\nat$ so that $x_n^*(x_i)>1-2\ep$ for $i\le n$. 
Now for $n\ge2$ of course we have $x_n^* \in C^X$. 
Suppose that $x_n^* = {1\over f(\ell_n)} \sum_1^{\ell_n} b_j^* \in 
C_{\ell(n)}^X$. 
If $\ell(n) \ge2$ then since $f(\ell_n) >1$ it must be true that $b_j^*(x_i) 
\ne0$ for at least two $j$'s for each $i\le n$ (provided $\ep>0$ 
satisfies ${1\over f(2)} <1-2\ep$). 
Thus if $\ell (n)\ge2$ for all $n$ we have that $\ell (n)\to\infty$. 
But then $x_n^* (x_1)\to0$ as $n\to\infty$, a contradiction. 
We are left with the case (passing to a subsequence) where $\ell(n)=1$ for 
all $n$ and so $x_n^* \in C_1^X = B^X$. 
But then $x_n^* \in B_{j(n)}^X$ where $j(n)\ge2$ for all $n$ and so by 
the argument we just gave adapted with minor notational changes, we obtain 
$x_n^* (x_1)\to0$, a contradiction.\qed
\enddemo

Henceforth $X= (c_{00},\dotNorm)$ shall denote the normed space obtained 
in Proposition~3.1. 
We shall write $\tilde A_m$, $A_m$, $B_m$, $C_m$, $B$, $C$, $\tilde\gA$, 
$\gA$ and $\gB$ rather than $\tilde A_m^X,A_m^X,\ldots$ etc. 

We can now state the main result of this paper. 

\proclaim{Theorem 3.2} 
Let $(y_i)$ be any block basis of $(e_i)$ in $X$. 
Let $(f_i)_1^n$ be any finite monotone basis. 
Then $(f_i)_1^n$ is block finitely represented in $(y_i)$. 
\endproclaim 

First we reduce the proof to the consideration of a special class of 
monotone bases. 
For $n\in\nat$ let $(S^n(i,j))_{i,j=1}^n$ be linearly independent vectors 
in some linear space. 
For $(a(i,j))_{i,j=1}^n \subseteq \real$ define 
$$\Big\| \sum_{i,j=1}^n a(i,j) S^n (i,j)\Big\| 
= \max_{j\le n} \max_{k\le n} \Big| \sum_{i=1}^k a(i,j)\Big|$$ 
$(S^n(i,j))_{i,j=1}^n$ is a normalized monotone basis under the norm when 
ordered lexicographically: 
$S^n(1,1),S^n(1,2),\ldots,S^n(1,n),S^n(2,1),\ldots,S^n(n,n)$. 
For fixed $j\le n$, $(S^n(i,j))_{i=1}^n$ is 1-equivalent to the usual 
summing basis of length~$n$.

\proclaim{Proposition 3.3} 
If for all $n\in\nat$, $(S^n(i,j))_{i,j=1}^n$ is block finitely 
represented in a basic sequence $(y_i)$ then  every finite monotone basis is 
block finitely represented in $(y_i)$. 
\endproclaim  

\demo{Proof} 
Let $(z_i)_{i=1}^m$ be a finite monotone basis for a space $Z$. 
We may assume that $Z$ is a subspace of $\ell_\infty^n$ for some $n$. 
For $i\le m$ write $z_i = (z(i,j))_{j=1}^n \in \ell_\infty^n$, 
i.e., $z_i(j) = z(i,j)$. 
For $i\le m$ set $w_i = \sum_{j=1}^n z(i,j)S^n(i,j)$. 
We claim that for all $(b_i)_1^m\subseteq \real$, 
$$\Big\| \sum_{i=1}^m b_iw_i\Big\| 
= \Big\| \sum_{i=1}^m b_i z_i\Big\|_{\ell_\infty^n}\ .$$ 
Indeed since $(z_i)_1^m$ is monotone, 
$$\eqalignno{
\Big\| \sum_{i=1}^m b_i z_i\Big\|_{\ell_\infty^n} 
& = \max_{k\le m} \Big\| \sum_{i=1}^k b_i z_i\Big\|_{\ell_\infty}^n 
 = \max_{k\le m} \max_{j\le n} \Big| \sum_{i=1}^k b_i z(i,j)\Big| \cr 
& = \Big\| \sum_{i=1}^m \sum_{j=1}^n b_i z(i,j) S(i,j)\Big\|
 = \Big\| \sum_{i=1}^m b_i w_i\Big\|\ .&\qed \cr}$$
\enddemo 

Thus we are reduced to proving that for all $n$, $(S^n(i,j))_{i,j=1}^n$ is 
block finitely represented in every block basis of $(e_i)$. 
Next we state our Main Lemma and show how it yields the theorem. 
The Main Lemma will be proved in \S4. 

\proclaim{Main Lemma 3.4} 
Let $Y$ be a block subspace of $X$. 
Then there exists a constant $C= C(Y) \ge1$ so that for all $\ep>0$ 
there exists $m_0\in\nat$ with the following property $(*)$ 
\vskip1pt
\item{$(*)$} For any $k\in \nat$, $\delta>0$ and $m\in \Im (\sigma)$ with 
$m\ge m_0$ and ${f(m)\over m} <\frac{\delta}8$ there exists $y\in Y$ 
with the following properties: 
\vskip1pt
\itemitem{\rm a)} $e_k<y$ 
\itemitem{\rm b)} $\|y\| \le C$
\itemitem{\rm c)} there exists $y^* \in B_m$ with $e_k^* <y^*$ and $y^*(y)>1$
\itemitem{\rm d)} {\rm i)} \ $\|y\|_{A_\ell} <\ep$ for all $\ell\in\nat$ 
with $\ell \ge m_0$
\itemitem{} {\rm ii)} $\|y\|_{B_\ell} < \ep$ for all $\ell\in\nat$ with 
$\ell \in[m_0,\bar m]\cup [\Bar m,\infty)$, 
\itemitem{} where $\bar m$ and $\Bar m$ are the predecessor of $m$ and 
successor of $m$ in $I_m(\sigma)$
\itemitem{\rm e)} {\rm i)}\ If $(x_1^* <\cdots < x_j^*)$ is 
$(\tilde\gA,M,\sigma)$-admissible with $\max\supp x_j^* \ge \min\supp y$ then 
$$|x^* (y)| <\delta\text{ for all } x^* \in \bigcup_{m\ge \sigma (x_1^*,
\ldots,x_j^*)} A_m$$ 
\itemitem{} {\rm ii)} If $(y_1^*,\ldots,y_j^*)$ is $(\gB,M,\sigma)$-admissible 
with $\max\supp y_j^* \ge \min\supp y$ then 
$$|x^* (y)| <\delta \text{ for all } y^* \in \bigcup_{m\ge \sigma 
(y_1^*,\ldots,y_j^*)} B_m$$ 
\endproclaim

Before proving Theorem 3.2 we present an elementary but often used lemma. 

\proclaim{Lemma 3.5} 
If $x_1<\cdots < x_m$ is a sequence in $B(X)$, $\ell\in\nat$, $\ell\ge2$ 
then 
$$\eqalign{\Big\| \sum_{i=1}^m x_i\Big\|_{A_\ell} 
& \le {1\over\ell} \max \biggl\{ \sum_{j=1}^\ell 
\Big\| \sum_{i=k_{j-1}}^{k_j} x_i \Big\| : 1\le k_0 \le k_1\le\cdots 
\le k_{\ell} \le m\biggr\}\cr 
& \le 1+{1\over\ell} \max \biggl\{ \sum_{j=1}^{\ell'} 
\Big\|\sum_{i=k_{j-1}+1}^{k_j} x_i\Big\| : 0\le k_0 <k_1<\cdots < k_{\ell'} 
\le m\text{ where } \ell'\le  \min (\ell,m)\biggr\}\cr}$$
\endproclaim 

\demo{Proof} 
Let $x^* = {1\over\ell} \sum_{j=1}^\ell x_j^*$ where $x_1^* <\cdots < 
x_\ell^*$ are in $B(X^*)\cap H$ and set $x= \sum_1^m x_i$. 
Let $k_0=1$ and for $j\ge1$ set 
$$k_j=\min \{i\ge 1 :\max\supp x_j^* \le \max\supp x_i\}$$ 
if such an $i$ exists and $k_j=m$ otherwise. 
Then $1=k_0\le k_1\le\cdots \le k_\ell \le m$ and 
$$\eqalignno{|x^*(x)| & = {1\over \ell} \Big| \sum_{j=1}^\ell x_j^*(x)\Big| 
= {1\over\ell} \Big| \sum_{j=1}^\ell \biggl( \sum_{i=k_{j-1}}^{k_j} 
x_j^* (x_i)\biggr)\Big|\cr 
& \le {1\over\ell} \sum_{j=1}^\ell \Big\| \sum_{i=k_{j-1}}^{k_j} x_i\Big\| \cr 
&\le {1\over\ell}\biggl[ \sum_{j=1}^\ell \|x_{k_{j-1}}\| 
+ \sum \Sb j=1\\ k_{j-1}<k_j\endSb ^\ell \Big\| \sum_{i=k_{j-1}+1}^{k_j} 
x_i\Big\|\biggr]\cr 
&\le 1+{1\over\ell} \max \biggl\{ \sum_{j=1}^{\ell'} \Big\| 
\sum_{i=k'_{j-1}+1}^{k'_j} x_i\Big\| : \ell' \le \min (\ell,m)\ ,\ 
0\le k'_0 <k'_1<\cdots < k'_{\ell'} \le m\biggr\}\ .&\qed\cr}$$
\enddemo

\proclaim{Corollary 3.6} 
Let $1\le p<\infty$, $k\in\nat$ and let $x\in X$ be an $\ell_p^k$-average 
with constant~2. Then for $\ell\ge2$, 
$$\|x\|_{A_\ell} \le k^{-1/p} + 2\ell^{-1/p}\ .$$ 
\endproclaim 

\demo{Proof} 
Let $x= k^{-1/p} \sum_{i=1}^k x_i$ where $x_1<\cdots <x_k$ is a 
normalized sequence 2-equivalent to the unit vector basis of $\ell_p^k$. 
By Lemma~3.5 there exist $0\le k_0<k_1<\cdots < k_{\ell'} \le k$ for 
some $\ell' \le \min (\ell,k)$ satisfying  
$$\eqalign{\|x\|_{A_\ell} 
& \le k^{-1/p} \left[ 1+{1\over\ell} \sum_{j=1}^{\ell'} 
\Big\| \sum_{i=k_{j-1}+1}^{k_j} x_i\Big\|\right]\cr 
&\le k^{-1/p} + \ell^{-1} k^{-1/p} \sum_{j=1}^{\ell'} 2(k_j-k_{j-1})^{1/p}\cr 
&\le k^{-1/p} + 2\ell^{-1} k^{-1/p} \ell' [k/\ell']^{1/p}\cr 
&\le k^{-1/p} + 2\ell^{-1} k^{-1/p} \ell^{1-1/p} k^{1/p}\cr 
&= k^{-1/p} + 2\ell^{-1/p}\ .\cr}$$ 
We have used the concavity of the function $x^{1/p}$ to obtain the third 
inequality.\qed
\enddemo 

\demo{Proof of Theorem 3.2} 
Fix $n$ and a block basis $(y_i)$ of $(e_i)$ in $X$. 
We shall prove that $(S^n(i,j))_{i,j=1}^n$ is block finitely represented 
in $(y_i)$. 
Let $\ep>0$ satisfy 
$$\gamma \equiv {1\over f(2)} \left[ 1+ {\ep\over 1-\ep}\right] <1
\text{ with } 2\ep <\gamma\ .$$ 
Let $C$ be as in Lemma 3.4 for $Y=[(y_i)]$. 
We next choose $k,\ell_0\in\nat$ and $\ep_0>0$ to satisfy 
$$\align
&{f(k)\over f(nk)} >1-\ep\ \text{ and }\ 6Cn {f(k)\over k} <\ep \ ,\tag 1\cr
&{4C n^3k^2\over f(\ell_0)} <\ep \ ,\ \text{ and}\tag 2\cr 
&6\ell_0 kn^2 \ep_0 <\ep\ .\tag 3\endalign$$ 

Then we let $m_0\ge M_{nk}$ be given by Lemma~3.4 for $\ep_0/2$. 
Let $P = \{(i,j,s) : 1\le i,j\le n$, $1\le s\le k\}$ be ordered 
lexicographically. 
Using $(*)$ in Lemma~3.4 we will recursively choose  for each $(i,j,s)\in P$ 
a block $y(i,j,s)\in \text{span}(y_t)$, $y^* (i,j,s)\in X^* \cap c_{00}$ 
and an interval $I(i,j,s)$ in $\nat$ and $m(i,j,s)\in\nat$ to satisfy 
conditions (4)--(11): 
\vskip1pt
\itemitem{(4)} $\supp y(i,j,s)\subseteq I(i,j,s)$\ ,
\itemitem{} $\supp y^* (i,j,s)\subseteq I(i,j,s)$\ ,
\vskip1pt
\noindent and $I(i,j,s) <I(i',j',s')$ if $(i,j,s) < (i',j',s')$. 
\vskip1pt
\itemitem{(5)} $\|y(i,j,s) \| \le C$\ ,
\itemitem{(6)} $y^* (i,j,s) (y(i,j,s)) =1$\ , 
\itemitem{(7)} 
$y^* (i,j,s) \in B_{m(i,j,s)}$ and for each $j\le n$ the 
family $(y^* (i,j,s))_{i\le n,\, s\le k}$ is $(\gB,M,\sigma)$-admissible 
(ordered lexicographically),  
\itemitem{(8)} 
a) $m(i,j,s)\ge m_0$ and if $(i,s) \ne (1,1)$ then 
$m(i,j,s) = \sigma (y^* (r,j,t): (r,j,t) <(i,j,s))$ and 
$${f(m(i,j,s)) \over m(i,j,s)} \le {1\over 8}\ {\ep_0 \over 
f^{-1}  \left( {1\over \ep_0} \max\supp y(i', j,s')\right) }$$
where $(i',s')$ is the predecessor of $(i,s)$,  
\itemitem{} b) 
$m(i,j,s)\ne m(r,\ell,t)$ if $(i,j,s)\ne (r,\ell,t)$ 
\itemitem{(9)} 
$\|y(i,j,s)\|_{A_\ell} <\ep_0$ if $\ell \ge m_0$\ ,  
\itemitem{(10)}
$\|y(i,j,s)\|_{B_\ell} <\ep_0$ if $\ell \in [m_0,\bar m]
\cup [\Bar m,\infty)$ where $\bar m$ and $\Bar m$ are the predecessor and 
successor of $m(i,j,s)$ in $\Im (\sigma)$. 
\itemitem{(11)} 
If $(i,s)\ne (1,1)$ and $(i',s')$ is the predecessor of $(i,s)$ then for any 
$(\gB,M,\sigma)$ or $(\tilde\gA,M,\sigma)$-admissible sequence 
$(x_1^* ,\ldots,x_{j_0}^*)$ with $\max\supp x_{j_0}^* \ge \min \supp 
y(i,j,s)$ then 
$$|x^* y(i,j,s)| < {\ep_0\over f^{-1}\left( {1\over\ep_0} \max\supp y
(i',j,s')\right)} \text{ for } x^*\in 
\bigcup_{t\ge\sigma (x_1^*,\ldots,x_{j_0}^*)} (B_t\cup A_t)\ .$$ 
\smallskip

\noindent 
Indeed let $(i,j,s)\in P$ and assume that $y(r,\ell,t)$, $y^*(r,\ell,t)$, 
$I(r,\ell,t)$ and $m(r,\ell,t)$ have been selected for all $(r,\ell,t) <
(i,j,s)$ so that (4)--(11) are satisfied for all such $(r,\ell,t)$ 
and furthermore if $(i,s) \ne (1,1)$ then defining 
$$m(i,j,s) \equiv \sigma (y^* (r,j,t): (r,t) <(i,s))$$ 
8) is satisfied as well for $(i,j,s)$. 

We then apply $(*)$ of Lemma 3.4 for the following parameters. 
If $(i,j,s)=(1,1,1)$ we let $k=k(1,1,1)=m_0$. 
Otherwise set 
$$k=k(i,j,s) = \max\supp y^* (i_0,j_0,s_0) \vee \max\supp y(i_0,j_0,s_0)$$ 
where $(i_0,j_0,s_0)$ is the predecessor of $(i,j,s)$ in $P$. 
If $(i,s) = (1,1)$ we take $\delta = \ep_0$ and $m_0\le m\equiv m(i,j,s) \in 
\Im (\sigma)$ to satisfy ${f(m)\over m} <{\delta\over 8}$ and to be 
distinct from all $m(r,\ell,t)$'s previously chosen. 
If $(i,s)\ne (1,1)$ we take 
$$\delta = {\ep_0\over  f^{-1}\left( {1\over \ep_0} \max\supp y(i',j,s')
\right)}$$ 
where $(i',s')$ is the predecessor of $(i,s)$ and $m\equiv m(i,j,s)\equiv 
\sigma (y^* (r,j,t) :(r,t) <(i,s))$. 
Note that in this case we have by our hypothesis (8) that 
${f(m)\over m} <{\delta\over8}$ and so $(*)$ of Lemma~3.4 does apply to 
these parameters: $(\ep,m_0,k,\delta,m) 
= (\ep_0,m_0,k(i,j,s),\delta,m(i,j,s))$. 

We thus obtain $\tilde y(i,j,s)>e_{k(i,j,s)}$ which satisfies 
$\|\tilde y(i,j,s)\|\le C$ and both (9), (10)  and (11) hold for 
$\tilde y(i,j,s)$ replacing $y(i,j,s)$. 

Furthermore there exists $y^* \in B_{m(i,j,s)}$ (and thus infinitely many 
such $y^*$'s) satisfying $y^* (\tilde y(i,j,s)) >1$ and 
$y^* > e_{k(i,j,s)}^*$. 
In particular we can choose $y^* (i,j,s)\in B_{m(i,j,s)}$ to be one of 
these $y^*$'s so that in addition if $(i,s)\ne (n,k)$ and $(i'',s'')$ is 
the successor of $(i,s)$ then $m(i'',j,s'') \equiv \sigma (y^* (r,j,t): 
(r,t) \le (i,s))$ also satisfies the condition in (8). 

We then set $y(i,j,s) = {\tilde y(i,j,s)\over y^*(\tilde y(i,j,s))}$ and 
$I(i,j,s)=(k(i,j,s),\max\supp y(i,j,s)\vee \max\supp y^* (i,j,s)]$. 
This completes the construction of the $y(i,j,s)$, $y^*(i,j,s)$, $I(i,j,s)$ 
and $m(i,j,s)$ satisfying (4)--(11). 

For $1\le i,j\le n$ define $y(i,j)={f(k)\over k} \sum_{s=1}^k y(i,j,s)$. 
>From (4) it follows that $\{y(i,j):1\le i,j\le n\}$ is a block basis of $(y_t)$. 
Let $(a(i,j))_{i,j\le n} \subseteq \real$ with $\|\sum_{i,j} a(i,j) 
S^n(i,j)\|=1$. 
We shall first show that for $y= \sum_{i,j} a(i,j)y(i,j)$, $\|y\|\ge1-\ep$. 
Fix $j_0,\ell\le n$ with $1=\sum_{i=1}^\ell a(i,j_0)$ (if the sum is $-1$ 
we replace all $a(i,j)$'s by $-(a(i,j))$. 
Define 
$$y^* = {1\over f(\ell k)} \sum_{i=1}^\ell \sum_{s=1}^k y^* (i,j_0,s)\ .$$
By (7) we have $y^* \in C_{\ell k}$ and furthermore
$$y^*(y) = {1\over f(\ell k)} {f(k)\over k} \sum_{i=1}^\ell \sum_{s=1}^k 
a(i,j_0)$$ 
by (4) and (6). But this 
$$= {1\over f(\ell k)} \ {f(k)\over k}\cdot k = {f(k)\over f(\ell k)} 
>1-\ep$$ 
by (1). 
It remains to prove that $\|y\| \le 1+\ep$.
\enddemo 

\proclaim{Claim 1} 
Let $x^* = {1\over f(r)} \sum_{t=1}^r x_t^* \in C_r$ with $r\ge 2$. Then 
$$|x^* (y)| \le 1+\ep \ \text{ or }\ |x^* (y)| \le \gamma \|y\|\ .$$ 
\endproclaim 

\noindent {\it Case 1\/}. $r\ge\ell_0$. 

For $j\le n$, the family $\{C^{-1} y(i,j,s) :i\le n$, $s\le k\}$ satisfies 
the conditions of Lemma~2.2a) for $(k,m,\ep) = (m_0,nk,\ep_0)$ and 
$(\D,M,\sigma) = (\gB,M,\sigma)$ or $(\tilde\gA,M,\sigma)$. 
Since $\|(a(i,j))\|_{\infty}\le 2$ we deduce from the second inequality in 
(2.2.2) that for $j\le n$. 
$$\eqalign{
\Big| \sum_{i=1}^n \sum_{s=1}^k a(i,j) x^* y(i,j,s)\Big| 
& \le {1\over f(r)} 2Cnk + {\min (r,kn)\over f(r)} 2Cnk\cr 
&\qquad + 2C\left[ 1+ {1\over f(r)} +2\ep_0 +1\right]\cr 
&\le {1\over f(r)} 4Cn^2 k^2 +6C\ .\cr}$$ 

It follows that since $r\ge\ell_0$ 
$$\eqalign{ |x^* (y)| & \le {f(k)\over k} \sum_{j=1}^n 
\Big| \sum_{i=1}^n \sum_{s=1}^k a(i,j) x^* y(i,j,s)\Big| \cr 
&\le {4Cn^3 k^2\over f(r)} + {6Cn f(k)\over k} <2\ep < \gamma\|y\|\ .\cr}$$
by (1) and (2)  and our choice of $\ep$ and $\gamma$. 
\medskip

\noindent {\it Case 2\/}. $r\le \ell_0$. 

For $t\le r$ let $x_t^* \in B_{m_t}$. 
It may be that $m_t = m(i,j,s)$ for some $t>1$ and $(i,j,s)\in P$ with 
$(i,s)\ne (1,1)$. 
In this case let $t_1$ be the maximum of such $t$'s and note that 
$$m_{t_1} = \sigma (x_1^*,\ldots,x_{t_1-1}^*) = m(i_1,j_1,s_1)$$ 
for some $(i_1,j_1,s_1)\in P$ with $(i_1,s_1) \ne (1,1)$. 
Also then 
$$m_{t_1} = \sigma (y^* (i,j_1,s) : (i,s) <(i_1,s_1))\ .$$ 
By the injectivity of $\sigma$ we deduce that  
$$(x_1^*,\ldots,x_{t_1-1}^*) = (y^* (1,j_1,1),y^* (1,j_1,2),\ldots, 
y^* (i_0,j_1,s_0)$$ 
where $(i_0,s_0)$ is the predecessor of $(i_1,s_1)$. 
Furthermore $m_t \ne m(i,j,s)$ for all $(i,j,s)\in P$ 
and $m_t >m_{t_1}\ge m_0$ if $t>t_1$. 
Thus by (10) we have that for all $(i,j,s)\in P$, 
$|x_t^* (y(i,j,s))| \le \ep_0$ if $t>t_1$ or if $t=t_1$ and 
$(i,j,s)\ne (i_1,j_1,s_1)$. 

>From these observations we obtain 
$$\eqalign{
|x^*(y)| & \le {1\over f(r)} \Big| \sum_{i=1}^{i_0-1} a(i,j_1) f(k) 
+ a(i_0,j_1) \sum_{s=1}^{s_0} {f(k)\over k}\Big|\cr 
&\qquad + {1\over f(r)} |x_{t_1}^* (y)| + {1\over f(r)}     
\sum_{t=t_1+1}^r |x_t^* (y)| \cr 
&\le \max \biggl\{ \Big| \sum_{i=1}^{i_0-1} a(i,j_1)\Big|,
\Big| \sum_{i=1}^{i_0} a(i,j_1)\Big| \cr 
&\qquad + {1\over f(r)} |x_{t_1}^* (y)| + {1\over f(r)} \ {f(k)\over k} 
\cdot 2\cdot n^2 k\ep_0 r\ .\cr}$$ 
The first term in the last inequality is obtained by noticing that if 
$i_0>1$ then necessarily $r\ge k$ while if $i_0=1$ then 
$${1\over f(r)} \sum_{s=1}^{s_0} {f(k)\over k} = {f(k)s_0\over f(r)k} 
\le {f(k)s_0\over f(s_0)k}$$ 
since $r\ge s_0$. 
The latter is not bigger than 1  since $s_0\le k$. 
Thus 
$$\eqalign{|x^*(y)| & \le 1+ {1\over f(r)} \ {2f(k)\over k} 
|x_{t_1}^* y(i_1,j_1,s_1)| \cr 
&\qquad + {1\over f(r)} \ {2f(k)\over k} n^2 k\ep_0
+ {2n^2 \ep_0r f(k)\over f(r)}\cr 
&< 1+ {1\over f(r)} \left[ {2Cf(k)\over k} + 2 n^2 f(k)\ep_0 +2n^2f(k)\ep_0
r\right]\cr 
&< 1+\ep\ \text{ using $r\le \ell_0$, (1) and (3).}\cr}$$
It remains to check the case in which for every $j\in \{1,\ldots,n\}$ 
and every $t>1$, $m_t \ne m(i,j,s)$ whenever $(i,s)\ne (1,1)$. 
In that case we obtain from (10) for $j\in \{1,\ldots,n\}$ that 
$|x_t^* y(i,j,s)| <\ep_0$ whenever $(i,s)\ne (1,1)$ and $t>t_0$ where $t_0$ 
is the smallest $t$ for which $\max\supp x_t^*\ge m_0$ (note that 
$m_{t_0+1} = \sigma (x_1^*,\ldots,x_{t_0}^*) > \max\supp x_{t_0}^*\ge m_0$). 

Thus we get 
$$\eqalign{ |x^*(y)| & = {1\over f(r)} \left| x_{t_0}^*(y) 
+ \sum_{t=t_0+1}^r x_t^* \biggl( \sum_{i,j=1}^n a(i,j) {f(k)\over k} 
\sum_{s=1}^k y(i,j,s)\biggr)\right| \cr 
& \le {\|y\|\over f(r)} + 2n^2 rf(k) {\ep_0\over f(r)}\cr 
& \le {\|y\|\over f(r)} + 2n^2 \ell_0 k{\ep_0\over f(r)}\cr 
& \le {1\over f(r)} (\|y\|+\ep)\ ,\quad\text{by 3)}\cr 
& \le \gamma \|y\|\ ,\cr}$$
by the choice of $\gamma$ and the fact that $\|y\|>1-\ep$. 
This completes Claim~1. 

\proclaim{Claim 2} 
Let $x^* = {1\over f(r)} \sum_{i=1}^r x_i^* \in B_r$ with $r\ge 2$. 
Then $|x^*(y)| \le \gamma \|y\|$. 
\endproclaim 

Indeed the case $r >\ell_0$ is handled exactly the same way as Case~1 in 
Claim~1. If $r\le \ell_0$ let $t_0\in\nat$ be minimal so that 
$\max\supp x_{t_0}^* \ge m_0$. 
Let $x_t^* \in A_{s_t}$. 
For $t>t_0$, $s_t >m_0$ by $(\sigma_2)$ and thus by (9), 
$|x_t^* y(i,j,s)| <\ep_0$ for $(i,j,s)\in P$. 
Since $e_{m_0} <y$ we deduce, using (3), that 
$$\eqalign{|x^* (y)| & \le {|x_{t_0}^* (y)| \over f(r)} + {1\over f(r)} 
\ell_0 n^2 k {f(k)\over k} \ep_0\cr 
&< {\|y\|+\ep \over f(r)} \le \gamma \|y\|\ .\cr}$$  
By the definition of $\dotNorm$ in $X$ we obtain $\|y\| \le 1+\ep$.\qed
\medskip

The proof of Theorem 3.2 yields the following corollary. 
Recall \cite{MMT} that if $Y$ has a basis $(y_i)$, $n\in\nat$ and 
$(x_i)_1^n$ is a normalized monotone basis then $(x_i)_1^n\in\{Y\}_n$ 
if $\forall\ \ep>0$ 
\vskip1pt
\itemitem{} $\forall\ k_1 \ \exists\ \ell_1 > k_1\ \exists\ z_1\in 
\text{span}(y_i)_{k_1}^{\ell_1}\ \forall\ k_2 >\ell_1\ \exists\ \ell_2>k_2$ 
$\exists\ z_2\in\text{span}(y_i)_{k_2}^{\ell_2}\ldots\ \forall\ 
k_n>\ell_{n-1}\ \exists\ \ell_n >k_n$ 
$\exists\ z_n\in\text{span}(y_i)_{k_n}^{\ell_n}$  with 
$(z_i)_1^n$ $1+\ep$-equivalent to $(x_i)_1^n$. 

\proclaim{Corollary 3.7} 
For all block subspaces $Y$ of $X$ and for all $n$, $\{Y\}_n$ is the set 
of all normalized monotone bases of length~$n$. 
\endproclaim 

\head \S4. Proof of the Main Lemma\endhead 

Since the proof of Lemma 3.4 is quite technical we first outline the argument. 
Let $Y$ be an arbitrary block subspace of $X$. 

\demo{Step 1} 
We first show that for some $1\le p<\infty$, $\ell_p$ is block finitely 
represented in $Y$. 
Indeed Krivine's theorem insures that there is a $p\in[1,\infty]$ so that 
$\ell_p$ is block finitely representable in $Y$. 
Secondly, we will observe (Lemma~4.1) that if $p=\infty$, then blocks of 
certain $\ell_\infty$-averages will produce for a given $k\in \nat$ and 
$\ep >0$ a sequence of length $k$ which is $(1+\ep)$-equivalent to the 
$\ell_1^k$-unit vector basis.
\enddemo 

\demo{Step 2} 
Let $1\le p<\infty$ be as found in Step~1. 
We first estimate the $\|\cdot\|_{B_\ell}$- and $\|\cdot\|_{A_\ell}$-norm of 
linear combination of certain $\ell_p$-averages (Lemma~4.3). 
Then we consider a sequence $(y_i)$, where $y_i$ is an $\ell_p^{k_i}$-average 
of constant $(1+\ep_i)$ with $k_i\uparrow \infty$ and 
$\ep_i\downarrow 0$. 
Let $E$ be a spreading model of a subsequence of $(y_i)$. 
\enddemo 

Either $c_0$ is block finitely representable in $E$. 
In that case we  will (Lemma~4.4) not only deduce that $c_0$ is block 
finitely representable in $[y_i]$ but also that we can choose for any 
$\ep>0$ and $k\in\nat$ an $\ell_\infty^k$ average $x$ of constant 
$(1+\ep)$, so that for any $(\tilde{\gA},M,\sigma)$ admissible sequence 
$(x_1^*,\ldots,x_j^*)$ with $\max\supp (x_j)>\min \supp x$ and any $x^*
\in \bigcup_{t\ge \delta (x_1^*,\ldots,x_j^*)} A_t$ we have $|x^*(y)|<
1+\ep$. 
This last condition says that $x$ is a ``good $\ell_\infty^k$-average'' 
but for $k'\gg k$ (where ``$k'\gg k$'' depends on $\min\supp x$) $x$ 
is a ``bad $\ell_\infty^{k'}$-average.'' 
We will call such a vector $x$ a special $\ell_\infty^k$-average of 
constant $1+\ep$. 

If $c_0$ is not block finitely representable in $E$ then for some 
$1\le q<\infty$, $\ell_q$ is block finitely representable in $E$. 
In this case we will be able to find a sequence $(z_k)$ in $y$ 
consisting of increasing $\ell_q$-averages. 
Furthermore $(z_k)$ satisfies the assumptions of Lemma~2.2(b) with 
$(\D,M,\delta)$ replaced by $(\B,M,\sigma)$ 
as well as by $(\tilde{\gA},M,\sigma )$. 
Applying Lemma~2.2(b) will give us that by replacing $p$ by $q$ and 
$(y_k)$ by $(z_k)$ we find ourselves in the first case. 

\demo{Step 3} 
Now we consider a spreading model of a sequence $(y_n)$ consisting 
of special $\ell_\infty^{k_n}$-averages of constant $(1+\ep_n)$, 
where $k_n\uparrow \infty$ and $\ep_n\downarrow 0$. 
Once again we have to distinguish between two cases. 
\enddemo 

\demo{Case 1} 
Up to passing to a subsequence we find a $C>0$ so that 
$$\varlimsup_{n_1\to\infty}\ \varlimsup_{n_2\to\infty}\ldots 
\varlimsup_{n_m\to\infty} {f(m)\over m} \Big\| \sum_{i=1}^m y_{n_i}\Big\| 
\le C\ ,$$ 
for all $m\in \nat$. 
\enddemo 

\demo{Case 2} 
Up to passing to a subsequence of $(y_n)$ we find $c_k\downarrow 0$ and $m_k
\uparrow \infty$ in $\nat$ so that 
$$\varlimsup_{n_1\to\infty} \ldots \varlimsup_{n_{m_k}\to\infty} c_k 
{f(m_k)\over m_k} \Big\| \sum_{s=1}^{m_k} y_{n_s} \Big\| =1\ .$$ 
\enddemo 

In the second case we let $(z_n)$ be a block sequence of the form 
$$z_n = c_n {f(m_n)\over m_n} \sum_{s=1}^{m_n} y_{k(n,s)}$$ 
and observe that $(z_n)$ satisfies the conditions of Lemma~2.2(b) for 
$\D = \tilde{\gA}$ as well as for $\D=\gB$, and deduce that $c_0$ is 
a spreading model of a subsequence of $(z_n)$. 
Taking  a block sequence $(\tilde y_n)$ of the form 
$$\tilde y_n = \sum_{i=1}^{k_n} z_{m(n,i)}$$ 
with $k_n\uparrow \infty$ we will observe that $(\tilde y_n)$ satisfies 
the hypothesis of Case~1. 
Thus we can assume Case~1 to be satisfied.

In that case we will show that choosing $C(\bar y)=C$ and letting 
$\ep>0$ and taking $m_0$ sufficiently large we can choose for any 
$k\in\nat$, $\delta>0$ and $m\in \Im (\delta)$, with $m\ge m_0$ and 
${f(m)\over m} < {\delta\over8}$, a vector $y\in Y$ to be of the form 
$$y= {f(m)\over m} \sum_{i=1}^m y_{n_i}\ ,$$ 
in order to satisfy the claim of Lemma~3.4.

We begin with an easy but important result. 

\proclaim{Lemma 4.1} 
Let $(y_i)$ be a block basis of $(e_i)$ in $X$. 
If $c_0$ is block finitely represented in $(y_i)$ then so is $\ell_1$. 
\endproclaim 

\demo{Proof} 
Given $n$ fixed we may choose a normalized block basis $(x_i)_{i=1}^n$ 
of $(y_i)$ along with functionals $x_i^*\in A_{m_i}$ so that 
\roster
\item"i)" $x_1^* <x_2^* <\cdots < x_n^*$; $x_i^* (x_i)>1/3$ and 
$x_i^* (x_j) =0$ for $1\le i\ne j\le n$
\item"ii)" For all $1\le k_1<k_2<\cdots <k_\ell \le n$ and all choices 
of sign $\pm$. 
$$(\pm x_{k_1}^* ,\pm x_{k_2}^* ,\ldots, \pm x_{k_\ell}^*)\text{ is 
$(\tilde\gA,M,\sigma)$-admissible.}$$
Indeed each $x_i$ will be a $(1+\ep)-\ell_\infty^{m_i}$ normalized average 
for suitable $m_i$ and small $\ep$. 
Thus if $(a_i)_1^\ell \subseteq \real$ and $1\le k_1<\cdots <k_\ell \le n$ 
and $\ep_i = \text{sign }a_i$ we have $x^* = {1\over f(\ell)} 
\sum_{i=1}^\ell \ep_i x_{k_i}^* \in B_\ell$ and 
$$\Big\| \sum_1^\ell a_i x_{k_i}\Big\| \ge x^* \biggl( \sum_{i=1}^\ell 
a_i x_{k_i}\biggr) \ge {1\over 3f(\ell)} \sum_{i=1}^\ell |a_i|\ .$$ 
\endroster

>From James' proof that $\ell_1$ is not distortable we obtain that $\ell_1$ 
is block finitely represented in $(y_i)$ \cite{J}.\qed 
\enddemo 

>From Lemma 4.1 and Krivine's theorem (\cite{K},\cite{L}) we have 

\proclaim{Corollary 4.2} 
If $(y_i)$ is a block basis of $(e_i)$ then there exists $p\in[1,\infty)$ 
so that $\ell_p$ is block finitely representable in $(y_i)$. 
\endproclaim

\proclaim{Lemma 4.3} 
Let $1\le p<\infty$, $0<\ep<1$ and $\ell\in\nat$. 
Let $(y_i)$ be a block basis of $(e_i)$ and let $k_1<\cdots <k_\ell$ 
satisfy 
$$\align
&\text{For $1\le i\le\ell$, $y_i$ is an $\ell_p^{k_i}$-average with 
constant $1+\ep$.}   \tag 4.3.1\cr 
&f\left( {k_1^{1/p^2} \ep^2\over 10\ell}\right) \ge {2\ell (1+2C_p)\over\ep}
\tag 4.3.2\endalign$$ 
and for $i=1,2,\ldots,\ell-1$, 
$${\ep\over2} f(k_{i+1}^{1/p}) \ge \sum_{s=1}^i |\supp y_i|\ .$$ 
Then for $(\alpha_i)_1^\ell \subseteq [-1,1]^\ell$, 
$y= \sum_{i=1}^\ell \alpha_i y_i$ and $m\ge2$, 
$$\text{if }\ x^* = {1\over f(m)} \sum_{j=1}^m x_j^* \in B_m\ ,
\leqno(4.3.3)$$ 
$$\align 
|x^*(y)| &\le \cases 
\ds {1\over f(m)} \Bigl( \max_{j\le m} |x_j^* (y)|+\ep)\ ,
&\text{$\ds m\le {\ep^2 k_1^{1/p^2} \over 10\ell}$}\cr 
\ds \max_{i\le \ell} |\alpha_i| (1+\ep) +\ep\ ,
&\text{$\ds m> {\ep^2 k_1^{1/p^2}\over 10\ell}$} \endcases \cr
&\le \cases 
\ds {1\over f(m)} (\|y\|+\ep)\ ,
&\text{$\ds m\le {\ep^2 k_1^{1/p^2}\over 10\ell}$}\cr
\ds \max_{i\le \ell} |\alpha_i| (1+\ep)+\ep\ ,
	&\text{$\ds m > {\ep^2 k_1^{1/p^2} \over 10\ell}$\ .}\endcases  
\endalign$$
{\rm (4.3.4)} For $x^* = {1\over m} \sum_{i=1}^m x_i^* \in A_m$, 
$$\eqalign{|x^*(y)| & \le \max_{s\le \ell} |\alpha_s| 
\left[ k_1^{-1/p} +4 \left( {\min (\ell,m)\over m}\right)^{1/p}\right]\cr 
&\qquad + {1\over m} \max \biggl\{ \sum_{j=1}^{\ell'} \Big\| 
\sum_{i\in E_j} \alpha_i y_i\Big\| :\ell' \le \min (\ell,m) 
E_1<\cdots <E_{\ell'} \text{ are intervals in } 
\{1,\ldots,\ell\}\biggr\}\ .\cr}$$
\endproclaim

\demo{Proof} 
Let $x^* = {1\over f(m)} \sum_{i=1}^m x_i^*\in B_m$ where $x_i^* = 
{1\over m_i} \sum_{j=1}^{m_i} x^* (i,j) \in A_{m_i}$. 
\enddemo 

\demo{Case 1} 
$m\le \ep^2 k_1^{1/p^2}/10\ell$. 

Let $j_0\in \{1,2,\ldots,m+1\}$ be maximal so that $\sum_{i=1}^{j_0-1} 
|\supp x_j^*| \le \ep k_1^{1/p}/2$. 
Thus if $j_0<m$ from $(\sigma_2)$ we have 
$m_{j_0+1} > \sum_{i=1}^{j_0} |\supp x_j^*| > {\ep k_1^{1/p}\over2}$. 
Thus by Corollary~3.6, 
$$\eqalign{ \sum_{j=j_0+1}^m |x_j^* (y)| 
&\le \sum_{i=1}^\ell \sum_{j=j_0+1}^m |x_j^* (y_i)| \cr 
&\le \ell m \left[ k_1^{-1/p} + 2m_{j_0+1}^{-1/p}\right]\cr 
&\le\ell m\left[ k_1^{-1/p} +2\cdot 2^{1/p} \ep^{-1/p} k_1^{-1/p^2}\right]\cr
&\le 5\ell m\ep^{-1} k_1^{-1/p^2} <\ep/2\ .\cr}$$ 
Also 
$$\eqalign{
\sum_{j=1}^{j_0-1} |x_j^*(y)| 
& \le \sum_{j=1}^{j_0-1} |\supp x_j^*|\ \|y\|_\infty\cr
&\le k_1^{-1/p} \sum_{j=1}^{j_0-1} |\supp x_j^*| <\ep/2\cr}$$ 
by our choice of $j_0$ and the fact that 
$$\|y\|_\infty  \le \max_i \|y\|_\infty \le \max_i k_i^{-1/p} = k_1^{-1/p}\ .$$
Thus 
$$|x^*(y)| \le {1\over f(m)} |x_{j_0}^* (y)| +\ep\ .$$ 
\enddemo 

\demo{Case 2} 
$m>\ep^2 k_1^{1/p^2} /10\ell$.

Choose $i_0\in \{1,\ldots,\ell+1\}$ maximal so that 
$$\sum_{i=1}^{i_0-1} |\supp y_i| < {\ep f(m)\over2}\ .$$ 
Then 
$$\Big| x^* \biggl( \sum_{i=1}^{i_0-1} \alpha_i y_i\biggr)\Big| 
< {\ep\over2}\ .$$ 
Also by (4.3.2) if $i_0 <\ell$, 
$${\ep f(m)\over2} \le \sum_{i=1}^{i_0} |\supp y_i| \le {\ep\over2} 
f(k_{i_0+1}^{1/p})$$ 
which yields $m\le k_{i_0+1}^{1/p}$. 
If $i_0 <i\le \ell$ we have 
$$|x^* (y_i)| \le {1\over f(m)} \sum_{j=1}^m |x_j^* (y_i)| 
\le {1\over f(m)} \sum_{j=1}^m (k_i^{-1/p} +2m_j^{-1/p})$$ 
by Corollary 3.6 and in turn by $(\sigma_3)$ this is 
$$\eqalign{
&\le {1\over f(m)} \left[ mk_{i_0+1}^{-1/p} +2C_p\right] 
\le {1\over f(m)} [1+2C_p]\cr 
&\le {1\over f(\ep^2 k_1^{1/p^2}/10\ell)} [1+2C_p] \le {\ep\over2\ell}\cr}$$ 
where the last inequality follows from (4.3.2). 

Thus $|x^* (\sum_{i=i_0+1}^\ell \alpha_i y_i)| \le \ep/2$. 
We obtain $|x^* (y)| \le \ep + |x^* (\alpha_{i_0} y_{i_0})|$ which 
completes the proof of (4.3.3), since $\|y_{i_0}\| < 1+\ep$.  

Let $x^* = \frac1m \sum_1^m x_i^* \in A_m$ and let $1=n_1<n_2<\cdots 
<n_{\ell+1}$ so that $\supp y_i \subseteq [n_i,n_{i+1})$ for $1\le i\le \ell$ 
and $n_{\ell+1} > \max\supp x_m^*$.  
For $1\le i\le\ell$ define $I_i = \{j :\supp x_j^* \subseteq [n_i,n_{i+1})\}$ 
and $m_i = |I_i|$. 
Note that $\sum_{i=1}^\ell {m_i\over m}\le1$. 
If $I_i\ne\emptyset$ then 
${1\over m_i} \sum_{j\in I_i} x_j^* \in A_{m_j}$ by $(\sigma_1)$ and so by 
Corollary~3.6, 
$$\eqalign{\Big| \sum_{j\in I_i} x_j^* (\alpha_i y_i)\Big| 
& = m_i {1\over m_i} \Big|\sum_{j\in I_i} x_j^* (y_i)\Big|\ |\alpha_i|\cr 
&\le m_i [k_i^{-1/p} + 2m_i^{-1/p}] |\alpha_i|\ .\cr}$$
Hence if ``$\mathop{{\sum}'}\limits_{i=1}^\ell$'' 
denotes ``$\sum\limits\Sb i=1\\ I_i\ne \emptyset\endSb^\ell$'' then 
$$\eqalign{ 
{1\over m} \mathop{{\sum}'}_{i=1}^\ell \Big| \sum_{j\in I_i} 
x_j^* (\alpha_i y_i)\Big| 
& \le \mathop{{\sum}'}_{i=1}^\ell {m_i\over m} k_i^{-1/p} |\alpha_i| 
+ 2\sum_{i=1}^{\ell'} {m_i\over m} m_i^{-1/p} |\alpha_i|\cr 
& \le \max_{s\le \ell} |\alpha_s| 
\left[ k_1^{-1/p} + 2\mathop{{\sum}'}_{i=1}^\ell {m_i\over m} 
m_i^{-1/p} \right]\ .\cr}$$ 
If $m\le \ell$ then we use the estimate, 
$$\mathop{{\sum}'}_{i=1}^\ell {m_i\over m} m_i^{-1/p} 
\le \max\Sb i\le \ell\\ I_i\ne\emptyset\endSb m_i^{-1/p} \le 1 \ .$$ 
If $m> \ell$ then by H\"olders inequality, 
$$\eqalign{ {1\over m} \mathop{{\sum}'}_{i=1}^\ell m_i^{1-1/p} 
&\le {1\over m} \biggl( \mathop{{\sum}'}_{i=1}^\ell 1^p\biggr)^{1/p} 
\biggl( \mathop{{\sum}'}_{i=1}^\ell m_i\biggr)^{1-1/p}\cr 
&\le {\ell^{1/p}\over m} m^{1-1/p} = \left( {\ell\over m}\right)^{1/p}\ .\cr}$$
Thus 
$${1\over m} \mathop{{\sum}'}_{i=1}^\ell \Big|\sum_{j\in I_i} 
x_j^* (\alpha_i y_i)\Big| 
\le \max_{s\le \ell} |\alpha_s| 
\left[ k_1^{-1/p} + 2\left( {\min (\ell,m)\over m}\right)^{1/p}\right] \ .$$
Let $I_0 = \{1,2,\ldots,m\}\setminus \bigcup_{i=1}^\ell I_i$ and 
$\ell' = |I_0|\le \min (\ell,m)$. 
Then for an appropriate choice of $k_1<\cdots < k_{\ell'}$ and intervals 
$E_1<\cdots <E_{m'} \subseteq \{1,\ldots,\ell\}$, 
$$\sum_{j\in I_0} |x_j^* (y)| 
\le \sum_{j=1}^{\ell'} |\alpha_{k_j}|\, \|y_{k_j}\| + 
\sum_{j=1}^{\ell'}    \Big\|\sum_{i\in E_j} \alpha_i y_i\Big\|\ .$$ 
Since 
$$\sum_{j=1}^{\ell'} |\alpha_{k_j}|\, \|y_{k_j}\| 
\le \min (\ell,m) \max_{s\le \ell} |\alpha_s| (1+\ep)\ ,$$ 
(4.3.4) follows from these estimates using that 
$$|x^*(y)| \le {1\over m} \sum_{i=1}^{\ell'} 
\Big| \sum_{j\in I_i} x^* (\alpha_i y_i) \Big| 
+ {1\over m} \sum_{j\in I_0} |x_j^* (y)|\ .\eqno\qed$$
\enddemo

\remark{Remark 4.4} 
By Corollary 4.2 for every block basis $(x_i)$ of $(e_i)$ there exists 
$1\le p<\infty$ so that for all $\ep>0$ and $\ell\in\nat$ there exists 
a block basis $(y_i)_{i=1}^\ell$ of $(x_i)$ satisfying (4.3.1) and (4.3.2). 
\endremark

\proclaim{Lemma 4.5} 
Let $(x_i)$ be a block basis of $(e_i)$, $\ep>0$ and $k\in\nat$. 
There exists $x\in\text{\rm span}(x_i)$ so that 
\vskip1pt
\noindent {\rm (4.5.1)} $x$ is an $\ell_\infty^k$-average with 
constant $1+\ep$ and 
\vskip1pt
\noindent {\rm (4.5.2)} if $(x_1^*,\ldots,x_j^*)$ is 
$(\tilde\gA,M,\sigma)$-admissible with $\max\supp x_j^* \ge\min\supp x$ then 
$$|x^*(x)|<\ep\text{ for all } x^* \in \bigcup_{t\ge \sigma (x_1^*,\ldots,
x_j^*)} A_t\ .$$
\endproclaim 

\demo{Proof} 
As in Remark 4.5 there exists $1\le p<\infty$ and a block basis $(y_i)$ 
of $(x_i)$ and a subsequence $(k_i)$ of $\nat$ satisfying for $\ep_i 
\equiv \ep/2^i$, 
\vskip1pt
\noindent {\rm (4.5.3)} For every $i$, $y_i$ is an $\ell_p^{k_i}$-average 
with constant $1+\ep_i$
\vskip1pt
\noindent {\rm (4.5.4)} For $\ell\in\nat$ and $\ell\le n_1<\cdots <n_\ell$
\vskip1pt
\itemitem{\rm a)} $f(k_{n_1}^{1/p^2} \ep_\ell^2 /10\ell) >2\ell 
(1+2C_p)/\ep_\ell$
\itemitem{\rm b)} ${\ep_\ell\over2} f(k_{n_{i+1}}^{1/p}) \ge 
\sum_{s=1}^i |\supp y_{n_s}|$ for $1\le i<\ell$ 
\itemitem{\rm c)} $(y_i)$ has a spreading model $(\tilde y_i)$ satisfying 
for $(\alpha_i)_1^\ell \in [-1,1]^\ell$, 
$$\bigg| \ \Big\| \sum_{i=1}^\ell \alpha_i y_{n_i}\Big\| 
- \Big\| \sum_{i=1}^\ell \alpha_i \tilde y_i\Big\|\ \bigg| <\ep_\ell\ .$$ 
\enddemo

Indeed we first choose a sequence $(y_i)$ satisfying a) and b) for all 
subsequences and then pass to a subsequence satisfying c). 

\demo{Case 1} 
$c_0$ is block finitely representable in $(\tilde y_i)$. 

Using c) we can thus find $N$ so that if $N<n_1<\cdots <n_N$ then there 
exists a normalized block basis $(w_i)_1^k$ of $(y_{n_i})_1^N$ which is 
$(1+\ep)$-equivalent to the unit vector basis of $\ell_\infty^k$. 
Thus $x=\sum_1^k w_i$ is an $\ell_\infty^k$-average with constant $1+\ep$. 
Now we do this choosing $n_1$ so large that $k_{n_1}^{-1/p}< \ep/6$ and 
if $(x_1^*,\ldots,x_j^*)$ is $(\tilde\gA,M,\sigma)$-admissible with 
$\max\supp x_j^* \ge \min\supp y_{n_1}$ then  
$\sigma (x_1^*,\ldots,x_j^*) > (24/\ep)^p N$. 

We can write $x=\sum_1^N \alpha_i y_{n_i}$ for some $(\alpha_i)_1^N
\subseteq [-2,2]^\nat$. 
Let $m\ge (24/\ep)^p N$ and let $x^*\in A_m$. 
>From (4.3.4) we obtain (we may assume $\ep/24 <1$),  
$$|x^*(x)| \le 2\left[k_{n_1}^{-1/p} +4\Bigl({N\over m}\Bigr)^{1/p}\right] 
+ {2N\over m} < 2\Big[ {\ep\over6} + {\ep\over6}\Big] + {2N\over m} <\ep$$ 
If Case  1 fails to hold then by Krivine's theorem (\cite{K}, \cite{L}) 
we have
\enddemo 

\demo{Case 2} 
$\ell_q$ is block finitely represented in $(\tilde y_i)$ for some 
$1\le q<\infty$. 

In this case we produce in $(y_i)$ a block basis $(z_i)$ of 
$\ell_q^{k'_i}$-averages with constant $1+\ep_i$ satisfying $k'_i\uparrow 
\infty$ and for all $\ell\le n_1<\cdots <n_\ell$

\itemitem{a)$'$} $f(k'_{n_1} {}^{1/q^2} \ep_\ell^2/10\ell) > 2\ell 
(1+2C_q)/\ep_\ell$ 

\itemitem{b)$'$} ${\ep_\ell\over2} f(k'_{n_{i+1}}{}^{1/q}) \ge 
\sum_{s=1}^i |\supp z_{n_s}|$ for $1\le i<\ell$ where $\supp (z_{n_s})$ 
is w.r.t. $(y_t)$ and 
$$z_i = \sum_{j=1}^{N_i} \alpha (i,j) y_{n(i,j)}$$ 
for some $N_i <n(i,1)<\cdots < n(i,N_i)$. 
The latter yields by (4.3.3) that for $i,m\in\nat$ and $x^* \in B_m$, 
$$|x^*(z_i)| \le {1\over f(m)} (\|z_i\|+\ep_i) + \ep_i 
+ (1+\ep_i) \max \{|\alpha (i,j)|: 1\le j\le N_i\}\ .$$
Since $z_i$ is an $\ell_q^{k'_i}$-average in $(y_s)$ with $k'_i\to\infty$ 
and $q<\infty$ it follows that 
$$\lim_{n\to\infty}\sup \{|x^* (z_i)| :  i,m\ge n\ ,\ x^*\in B_m\}=0\ .$$ 

Thus the hypothesis of Lemma 2.2b) is satisfied with $\D,M,\sigma$ 
replaced by $\gB,M,\sigma$ and $y_i$ replaced by ${z_i\over \|z_i\|}\to1$. 
Hence for all $(\alpha_i)_1^\ell\subseteq \real$, 
$$ \varlimsup_{n_1\to\infty} \ldots \varlimsup_{n_\ell\to\infty} 
\Big\| \sum_{i=1}^\ell \alpha_i z_{n_i}\Big\| 
\le  \max \biggl\{ \|(\alpha_i)\|_\infty, 
\max_{j\le \ell} \varlimsup_{n_j\to\infty} \ldots  \varlimsup_{n_\ell\to
\infty} \Big\| \sum_{i=j}^\ell \alpha_i z_{n_i}\big\|_B\biggr\}$$ 
But by (4.3.3), which applies due to a$'$) and  b$'$), 
$$ \varlimsup_{n_1\to\infty} \ldots \varlimsup_{n_\ell\to\infty} 
\Big\| \sum_{i=1}^\ell \alpha_i z_{n_i}\Big\|_B 
\le \max \biggl\{ \|(\alpha_i)\|_\infty , {1\over f(2)}  
\varlimsup_{n_1\to\infty}\ldots \varlimsup_{n_\ell\to\infty} 
\Big\| \sum_{i=1}^\ell \alpha_i z_{n_i}\Big\|\biggr\}$$ 
which together with the above inequality implies 
$$\varlimsup_{n_1\to\infty} \ldots \varlimsup_{n_\ell\to\infty} 
\Big\| \sum_{i=1}^\ell \alpha_i z_{n_i}\Big\| = 
\|(\alpha_i)\|_\infty\ .$$ 

The lemma follows by this and (4.3.4) if we set 
$x= \sum_1^k z_{n_i}/\|z_{n_i}\|$ for a suitable choice of 
$n_1<\cdots < n_k$.\qed
\enddemo 

\subhead Proof of the Main Lemma 3.4\endsubhead

By virtue of Lemma 4.5 we can choose a block sequence $(y_i)$ in $Y$ along 
with sequences $\ep_i\downarrow 0$ with $\ep_1<1/2$ and a subsequence 
$(k_i)$ of $\nat$ so that conditions (1) and (2) hold for all $i\in\nat$. 
\vskip1pt
\itemitem{(1)} a) $y_i$ is an $\ell_\infty^{k_i}$-average with constant 
$(1+\ep_i)$ 
\itemitem{} b) if $(x_1^*,\ldots,x_j^*)$ is $(\tilde \gA,M,\sigma)$-admissible 
with $\max\supp x_j^* \ge \min\supp y_i$ then 
$$|x^* (y_i)| <\ep_i\ \text{ for all }\ x^*\in \bigcup_{m\ge\sigma (x_1^*,
\ldots,x_j^*)} A_m$$ 
\itemitem{} c) there exists $(y_i^*) \subseteq B(X^*)\cap A$ with 
$\supp y_i^* \subseteq [\min\supp y_i,\max\supp y_i]$, 
$y_i^* (y_i) >1/3$ and $y_{i+1}^* \in A_{k_{i+1}}$ where $k_{i+1} \ge 
\sigma (y_1^*,\ldots, y_i^*)$ 
\itemitem{(2)}  ${1\over\ep_i} \max\supp y_i< f\left( {\ep_i\over 
\ep_{i+1}}\right)$ 

\noindent 
Note that any subsequence of $(y_i,\ep_i,k_i)_\nat$ also satisfies 
conditions (1) and (2) (for condition c) this uses $(\sigma_1)$). 

Let $m\in\nat$, $m<n_1<\cdots <n_m$, set $x_i=y_{n_i}$ for $i\le m$ and 
$x= \sum_1^m x_i$. 

We first obtain estimates for $\|x\|_{A_\ell}$ and $\|x\|_{B_\ell}$ for 
$\ell\ge2$. 
>From Lemma~3.5 and the fact that $\|x_i\|\le 1+\ep_m$ for $i\le m$ we have 
$$\|x\|_{A_\ell} \le (1+\ep)\left( 1+{m\over\ell}\right) < 2{\ell+m\over\ell}
\ .\tag 3$$

>From (1) b) and (2) and the fact that $\ep_i\downarrow$ 
it is easy to check that Lemma~2.2~a) applies for 
$\ep=\ep_m$ and $\D = \tilde\gA$, $k=m$ (and $y_i$ in Lemma~2.2 replaced 
by $x_i\over 1+\ep_m$). 
We obtain from the second estimate in (2.2.2) that 
$$\eqalign{\|x\|_{B_\ell} 
&\le \sup_{\ell'\ge \ell} {\|x\|_{A_{\ell'}} \over f(\ell)} 
+ {\min (\ell,m)\over f(\ell)} \sup_{\ell'\ge m} \|x\|_{A_{\ell'}} \cr 
&\qquad + (1+\ep_m)
\left[ 1+ {1\over f(\ell)} + 2\ep_m + {2\over f(\ell)}\right] \ .\cr}$$ 
We have used that $x_1^*\in A_t$ for $t\ge M_\ell\ge \ell$ if 
$(x_1^*,\ldots,x_j^*)$ 
is $(\tilde \gA,M,\sigma)$-admissible to obtain the first term.

Thus since $\sup_{\ell'\ge m}\|x\|_{A_{\ell'}}<4$ by (4) we obtain using 
$\ep_m <\ep_1 <{1\over2}$
$$\|x\|_{B_\ell} \le \sup_{\ell'\ge \ell} {\|x\|_{A_{\ell'}}\over f(\ell)} 
+ 9 {\min (\ell,m)\over f(\ell)} + 3\ . 
\tag 4$$ 

Let $(\tilde y_i)$ be a spreading model of a subsequence of $(y_i)$. 
It may be that for some constant $C$ we have that $C(m)\equiv {f(m)\over m} 
\|\sum_{i=1}^m \tilde y_i\| <C$ for all $m$. 
If  so  then by passing to a subsequence of $(y_i)$ we may assume that we 
have 

\demo{Case 1} 
For all $m<n_1<\cdots <n_m$, 
$${f(m)\over m} \Big\| \sum_1^m y_{n_i}\Big\| \le C\ .$$ 
\enddemo

In the remaining case we have $\varlimsup_m C(m)=\infty$. 
Select a subsequence $C(m_n)\uparrow \infty$ with for all $n$, 
$C(i) < C(m_n)$ if $i<m_n$. 
Thus 
$${1\over C(m_n)}\ {f(m_n)\over m_n} \Big\| \sum_{j=1}^{m_n} \tilde y_j
\Big\| =1$$ 
and if $m'<m_n$ then 
$${1\over C(m_n)}\ {f(m')\over m'} \Big\| \sum_1^{m'} \tilde y_j 
\Big\| <1\ .$$
Set $c_n \equiv C(m_n)^{-1}$. 
Thus $c_n\downarrow 0$. 
Hence if Case~1 fails to hold we have 

\demo{Case 2} 
There exists a block basis $(z_n)$ of $(y_i)$ where 
$$z_n = c_n {f(m_n)\over m_n} \sum_{s=1}^{m_n} y_{k(n,s)}$$ 
for some $m_n<k(n,1) <\cdots < k(n,m_n)$ and $\|z_n\|$ is chosen 
so that $\big|\,\|z_n\|-1\big| <\ep_n$ 
and if $F\varsubsetneqq \{1,\ldots,m_n\}$ then 
$$\Big\| \sum_{s\in F} y_{k(n,s)}\Big\| < {|F| \over c_nf(|F|)}\ .
\tag 5$$ 
\enddemo 

We return now to Case 1 and complete the proof of Lemma~3.4 in this 
situation. 
Let $\ep>0$. 
Choose $m_0$ so that (as usual $\bar m$ and $\Bar m$ are the predecessor and 
successor of $m$ in $\Im (\sigma)$) 
$$\align
&{C\over f(m_0)} + {f(m_0)\over m_0} < {\ep\over 4}\ , \tag 6\cr 
&\sup_{\ell,m\ge m_0} C{f(m)\over \ell f\left( {m\over \min (\ell,m)}\right)} 
< {\ep\over2}\ , \text{ and}\tag 7\cr 
&\sup \left\{ {f(m)\ \min (\ell,m)\over m\ f(\ell)} :\quad  
\eqalign{&m\in \Im (\sigma)\ ,\ m\ge m_0\ \text{ and}\cr 
\noalign{\vskip-4pt}
	&\ell \in [m_0,\bar m] \cup [\Bar m,\infty)\cr} 
\right\} < {\ep \over 36}\ ,  \tag 8
\endalign$$ 
where $(\sigma_4)$ is used to get (8). 

To verify $(*)$ of Lemma~3.4 we let $k\in\nat$, $\delta>0$ and 
$m\in \Im (\sigma)$ with $m\ge m_0$ and ${f(m)\over m} <{\delta\over8}$. 
Choose $n_0>\max\{k, M_m\}$ so that if $(x_1^*,\ldots,x_j^*)$ is 
$(\tilde\gA,M,\sigma)$-admissible with $\max\supp x_j^* \ge\min\supp y_{n_0}$ 
then $\sigma (x_1^*,\ldots,x_j^*)\ge m$ while if $(x_1^*,\ldots,x_j^*)$ 
were $(\gB,M,\sigma)$-admissible with $\max\supp x_j^* \ge \min\supp y_{n_0}$ 
then 
$$f(\sigma (x_1^*,\ldots,x_j^*)) > {20f(m)\over\delta}
\ .\tag 9$$ 
Choose $n_0<n_1<\cdots < n_m$. 
Thus, by Case 1,  
$$\Big\| {f(m')\over m'} \sum_{s\in F} y_{n_s}\Big\| \le C \text{ for } 
F\subseteq \{1,\ldots,m\}\ ,\ |F|=m'\ .
\leqno(10)$$ 

Define $y= {f(m)\over m} \sum_{i=1}^m y_{n_i}$. 
We have $\|y\|\le C$ and $e_k<y$. 
Also by (1)~c), $(\sigma_1)$ 
and the fact that $n_0>M_m$ there exists $y^*\in B_m$ 
with $e_k^* <y^*$ and 
$$y^* (y) > {1\over f(m)} \ {f(m)\over m} \sum_{i=1}^m {1\over3}={1\over3}\ .$$
We have verified a), b) of Lemma 3.4 
and c) with constant $1/3$ rather than $1$. 
However this ``weaker result'' will formally imply the stated version.  
It remains to check conditions d) and e). 

Condition d) i) follows from Lemma 3.5 and the choice of $m_0$. 
Indeed for $\ell\ge m_0$ 
$$\eqalign{\|y\|_{A_\ell} 
& \le 2{f(m)\over m} + {f(m)\over m} \ {1\over \ell} \max 
\biggl\{ \sum_{j=1}^{\ell'} {Cm_i\over f(m_i)} : 
\ell'\le \min (\ell,m)\ ,\ m_1+\cdots+m_{\ell'} =m\biggr\}\cr 
&\le 2{f(m)\over m} + {f(m)\over\ell m} \max_{\ell'\le\min (\ell,m)} 
\left\{ {Cm\over f\left({m\over\ell'}\right)}\right\}\cr 
&\le 2{f(m)\over m} + {Cf(m)\over\ell} \ {1\over f\left({m\over \min(\ell,m)}
\right)} < \ep\ .\cr}$$ 
We have used the concavity of $x/f(x)$ along with (6), (7), (10), 
$\|y_i\| <2$ and the fact that $m\ge m_0$. 

Let $\ell\in [m_0,\bar m]\cup [\Bar m,\infty)$. 
>From (4) 
$$\eqalign{\|y\|_{B_\ell} & \le {f(m)\over mf(\ell)} \sup_{\ell'\ge\ell} 
\Big\| \sum_1^m y_{n_i}\Big\|_{A_{\ell'}} 
+ 9 {\min (\ell,m)\over f(\ell)}\ {f(m)\over m} + {3f(m)\over m} \cr 
& \le {C\over f(\ell)}  + 9{f(m)\min (\ell,m)\over mf(\ell)} 
+ {3f(m)\over m}\cr
&< \ep\cr}$$ 
(using (10), (6), (8) and $m\ge m_0$, $\ell\ge m_0$). 
Thus d)~ii) holds. 

If $(x_1^*,\ldots,x_j^*)$ is $(\tilde\gA,M,\sigma)$-admissible with 
$\max\supp_j^*\ge \min\supp y$ then our choice of $n_0$ implies that 
$\sigma (x_1^*,\ldots,x_j^*)\ge m$. 
Furthermore for $\ell\ge m$ from (4) we have 
$$\|y\|_{A_\ell} \le {f(m)\over m} 2{\ell+m\over\ell} \le 
4{f(m)\over m} <\delta\ .
\tag 11$$
Thus (e) i) of Lemma 3.4 holds. 

If $(y_1^*,\ldots,y_j^*)$ is $(\gB,M,\sigma)$-admissible with 
$\max\supp y_j^* \ge \min\supp y$ then 
for $\ell\ge \sigma (y_1^*,\ldots,y_j^*)$ it follows from (4), 
(9) and (11) that 
$$\|y\|_{B_\ell} \le {\delta\over f(\ell)} + 9{f(m)\over f(\ell)} 
+3{f(m)\over m} <\delta\ .$$
Thus e) ii) holds. 

Now let $(z_n)$ be as described in Case~2 above. 
For $\ell\ge2$ by Lemma~3.5  and (5), 
$$\align
(12)\qquad \|z_n\|_{A_\ell} & \le {2c_nf(m_n)\over m_n} \cr  
&\qquad + {c_nf(m_n)\over \ell m_n}
\max \biggl\{ \sum_{j=1}^{\ell'} {k_j\over c_nf(k_j)} : 
\ell' \le \min (\ell,m_n)\ ,\ k_1+\cdots + k_{\ell'} =m_n\biggr\}\cr 
&\le {2c_n f(m_n)\over m_n} + {f(m_n)\over \ell m_n} \ 
{m_n\over f\left( {m_n\over \min (\ell,m_n)}\right) }\cr 
&= {2c_n f(m_n)\over m_n} + {f(m_n)\over \ell f\left({m_n\over 
\min (\ell,m_n)}\right)} \to {1\over\ell}
\endalign$$
as $n\to\infty$. 

Furthermore from (12) we obtain, 
$$\lim_{i\to\infty} \sup \{\|z_n\|_{A_\ell} : n,\ell \ge i\}=0\ . 
\tag 13$$
For $\ell\ge2$ by (4) 
$$\|z_n\|_{B_\ell} \le (1+\ep_n) {1\over f(\ell)} + 9c_n {f(m_n)\over m_n} 
\ {\min (\ell,m_n)\over f(\ell)} + 3c_n {f(m_n)\over m_n}\ ,
\tag 14$$ 
where we have used 
$$\|z_n\|_{A_\ell} \le \|z_n\| \le 1+\ep_n\ .$$ 
It follows that we have, using $c_n\to0$,  
$$\lim_{i\to\infty} \sup \{ \|z_n\|_{B_\ell} : \ell,n\ge i\}=0\ .
\tag 15$$

Hence Lemma 2.2~b) applies for $\D$ replaced by either $\tilde\gA$ or $\gB$ 
(we do not have $\|z_n\|\le1$ but rather $\|z_n\|\to1$ which suffices). 
For all $k$ and $(\alpha_i)_1^k\subseteq \real$
$$F((\alpha_i)_1^k) \equiv 
\varlimsup_{n_1\to\infty} \ldots \varlimsup_{n_k\to\infty} 
\Big\| \sum_{i=1}^k \alpha_k z_{n_i} \Big\| 
= \|(\alpha_i)\|_\infty\ .
\tag 16$$ 

We prove this by induction on $k$. 
For $k=1$ the result is obvious. 
Assume that (16) holds for $k'<k$ with $k>1$.  
>From (2.2.5), applied  twice, 
$$\eqalign{F((\alpha_i)_1^k) 
& \le \max \left\{ \|(\alpha_i)\|_\infty, 
\varlimsup_{n_1\to\infty} \ldots \varlimsup_{n_k\to\infty} 
\Big\| \sum_{i=1}^k \alpha_i z_{n_i}\Big\|_B\right\}\cr 
&\le \max \left\{ \|(\alpha_i)\|_\infty, 
\varlimsup_{n_1\to\infty} \ldots \varlimsup_{n_k\to\infty} 
\Big\| \sum_{i=1}^k \alpha_i z_{n_i}\Big\|_A\right\}\ .\cr}$$ 
By (13) we see that there exists $k'<\infty$ so that 
$$\varlimsup_{n_1\to \infty}\ldots \varlimsup_{n_k\to\infty} 
\Big\| \sum_{i=1}^k \alpha_i z_{n_i}\Big\|_A = 
\varlimsup_{n_1\to\infty} \ldots \varlimsup_{n_k\to\infty} 
\Big\| \sum_{i=1}^k \alpha_i z_{n_i}\Big\|_{\bigcup_{\ell=2}^{k'} A_\ell}\ .$$ 
There exists $2\le \ell\le k'$ so that this 
$$= \varlimsup_{n_1\to\infty} \ldots \varlimsup_{n_k\to\infty}  
\Big\|\sum_{i=1}^k \alpha_i z_{n_i}\Big\|_{A_\ell}\ .$$ 
>From Lemma 3.5 this limit is not bigger than 
$$\varlimsup_{n_1\to\infty} \ldots \varlimsup_{n_k\to\infty}  
{1\over \ell} \sum_{j=1}^\ell \Big\|\sum_{i=k_{j-1}}^{k_j} \alpha_i z_{n_i} 
\Big\|$$ 
for some $1\le k_0 \le k_1\le\cdots \le k_\ell \le k$. 
Thus from the induction hypothesis this is in turn either 
$\le\|(\alpha_i)\|_\infty$ if $k_i\in (1,k)$ for some $i$ or otherwise 
$$\le {\ell-1\over\ell} \|(\alpha_i)\|_\infty + {1\over\ell} 
F((\alpha_i)_1^k)\ .$$ 
In the latter case 
${\ell-1\over\ell} F((\alpha_i)_1^k)\le {\ell-1\over\ell} 
\|(\alpha_i)\|_\infty$ and so we deduce that (16) holds (the upper 
$\infty$-estimate implies the lower one). 

Using (16) we can construct a block basis $(\bar y_n)$ of $(z_n)$ of the 
form $\bar y_n = \sum_{i=1}^{\bar k_n} z_{m(n,i)}$ for some 
$\bar k_n\uparrow \infty$ and $\bar k_n < m(n,1) <\cdots <m(k,\bar k_n)$ 
with $\bar y_n$ being (essentially) an $\ell_\infty^{\bar k_n}$-average. 
There is a slight difficulty in that $\|z_n\|\to1$ as opposed to 
$\|z_n\|=1$ but we shall ignore this trivial obstacle. 
We may presume that for some $\ep_i\downarrow 0$, 

\noindent (17) for $i\in\nat$

\itemitem{a)} $\bar y_i$ is an $\ell_\infty^{\bar k_i}$-average of $(z_n)$ 
with constant $1+\ep_i$. 

\itemitem{b)} If $(x_1^*,\ldots,x_j^*)$ is $(\tilde\gA,M,\sigma)$ or 
$(\gB,M,\sigma)$-admissible with $\max\supp x_j^*\ge \min\supp \bar y_i$ 
then 
$$|x^* (\bar y_i)| < \ep_i\text{ for }x^* \in \bigcup_{t\ge\sigma(x_1^*,
\ldots,x_j^*)} A_t\cup \bigcup_{t\ge \sigma (x_1^*,\ldots,x_j^*)} B_t\ .$$ 

\itemitem{c)} There exist $y_i^* \in A_{\bar k_i}$ with 
$$y_i^* (\bar y_i) > \frac13\text{ and } k_{i+1} \ge \sigma 
(y_1^*,\ldots,y_i^*)\ .$$

Part b) is achieved via (13) and  (15). 
Also (17~c) yields that  

\noindent (18) 
If $M_m <n_1 <\cdots <n_m$ then 
$$\Big\| \sum_{i=1}^m \bar y_{m_i}\Big\| 
\ge \Big\| \sum_{i=1}^m \bar y_{m_i}\Big\|_{B_m} 
> {1\over3} \ {m\over f(m)}\ .$$ 

We can also assume the following growth condition. 
$$\align
\hbox{\rm a)}&\qquad 1-\ep_1 - {1\over f(2)} >0\ ,\text{ and} \tag 19\cr
\hbox{\rm b)}&\qquad\text{ For }i\in\nat\ ,\quad 
\ep_{i+1} < {\ep_i \over if^{-1} \left( {i\max\supp\bar y_i\over \ep_i}
\right)}\ .
\endalign$$

Conditions (17)--(19) yield that conditions (1)--(4) hold for the sequence 
$(\bar y_i)$ replacing $(y_i)$. 
We shall now show that the sequence $(\bar y_i)$ satisfies for some 
$C$ that for all $m$ there exists $n_0$ so that if $n_0<n_1<\cdots <n_m$ 
then $\|\sum_1^m \bar y_{n_i}\| \le C{m\over f(m)}$. 
Thus we return to Case~1 and the proof will be complete. 

Choose $m_0$ so that 
$$\align
\hbox{\rm a)}\qquad&{1\over f(2)} + {6f(m_0)\over m_0} < 1-\ep_1\ , \tag 20\cr 
\hbox{\rm b)}\qquad& \left( 1+2\ep_{m_0} +{1\over f(2)}\right) 
(1+\ep_{m_0}) +\ep_{m_0} <2\ ,\text{ and}\cr 
\hbox{\rm c)}\qquad&{m_0\over f(m_0)} >3\ .
\endalign$$ 
Define 
$$C= \max \left\{ 2m_0, {10\over 1-\ep_1- {1\over f(2)} }\right\}\ .$$ 
Our claim is trivial for $m\le m_0$. 
Let $m>m_0$. 
Using (17~b) 
choose $x=\sum_1^m x_i$ where $x_i=\bar y_{n_i}$, $m<n_1<\cdots <n_m$ 
and where $n_1$ is so large that 
$$\|\bar y_n\|_{B_\ell} < {\ep_m \over m^2+2}
\tag 21$$ 
for all $n,\ell\ge n_1$. 
We first show that for $\ell\ge2$, 
$$\|x\|_{C_\ell} < \|x\|(1-\ep_1)\ .
\tag 22$$ 

Conditions (17 b) and (19 b) imply that Lemma~2.2~a) holds 
for $\left( {\bar y_{n_i}\over 1+\ep_m}\right)_1^m$ in the setting 
$\D = \gB$, $\ep = \ep_m$ and $k=n_1$. 
Thus for $\ell\ge2$, by (2.2.2), 
$$\eqalign{\|x\|_{C_\ell} 
& \le {1\over f(\ell)} \|x\|_B + {\min (\ell,m)\over f(\ell)} 
\sup_{\ell\ge n_1} \|x\|_{B_\ell} \cr 
&\qquad + (1+\ep_m) \left( {1\over f(\ell)} + 1+2\ep_m\right) 
+ {2\over f(\ell)} \max_{i\le m} 
\sup_{\ell\ge n_1} \|x_i\|_{B_\ell} \cr 
&\le {1\over f(\ell)} \|x\|_B + (1+\ep_m)
\left( {1\over f(\ell)} + 1+2\ep_m \right) 
+ {\min (\ell,m)m+2\over f(\ell)}  \max \Sb i\le m\\ \ell\ge n_1\endSb 
\|x_i\|_{B_\ell} \cr 
&\le {1\over f(\ell)} \|x\|_B + (1+\ep_m) \left( {1\over f(\ell)} +1+2\ep_m
\right) +\ep_m\ .\cr}$$ 
We used (21) to get the last estimate. 
By (20~b) this is $\le {1\over f(\ell)} \|x\|_B +2$. 
Thus by (18), 
$$\align \|x\|_{C_\ell} & \le {1\over f(\ell)} \|x\|_B 
+ 6{f(m)\over m} \|x\|_B  \tag 23 \cr 
&\le \left[ {1\over f(\ell)} + {6f(m)\over m}\right] \|x\| 
\le (1-\ep_1) \|x\|\endalign$$
by (20 a). 

Finally if $\ell\ge2$ and if $\|x\|_{B_\ell} > (1- \ep_1)\|x\|$ 
then by (4) 
$$\eqalign{\|x\|_{B_\ell} 
& \le {\|x\|\over f(\ell)} + 3+ 9{\min (\ell,m)\over f(\ell)} \cr 
&\le {\|x\|\over f(\ell)} + 3+ 9 {m\over f(m)}\cr 
&\le {\|x\|\over f(\ell)} + 10 {m\over f(m)}\cr}$$ 
(using $3\le {m_0\over f(m_0)} < {m\over f(m)}$ by (20~c)). 
Thus 
$$(1-\ep_1) \|x\| - {1\over f(\ell)} \|x\| < {10m\over f(m)}$$ 
and so by our choice of $C$, 
$$ \|x\| \le  \left( {10\over 1-\ep_1-{1\over f(\ell)}}
\right) {m\over f(m)} \le C{m\over f(m)}\ .$$
This, thankfully, concludes the proof.\qed

\enddocument